\documentclass[12pt]{amsart}

\usepackage{amsmath,amsthm,mathrsfs,amssymb}

\newcommand{\N}{\mathbb{N}}
\newcommand{\R}{\mathbb{R}}
\newcommand{\Z}{\mathbb{Z}}

\newcommand{\NT}{\widetilde{N}}
\newcommand{\supp}{\mathrm{supp}}
\newcommand{\essinf}{\mathrm{ess\hspace{1mm}inf}}
\newcommand{\set}[1]{\left\{#1\right\}}
\newcommand{\e}{\varepsilon}

\renewcommand{\P}{\mathbb{P}_{\theta}}
\newcommand{\E}{\mathbb{E}_{\theta}}

\def\ds{\displaystyle}

 \newtheorem{thm}{Theorem}[section]
 
 \newtheorem{lem}[thm]{Lemma}
 \newtheorem{prop}[thm]{Proposition}
 \theoremstyle{definition}
 
 \theoremstyle{remark}
 \newtheorem{rem}[thm]{Remark}
 
 \numberwithin{equation}{section}

\begin{document}

\title[IDS for perturbed lattice]
{Classical and quantum behavior of the integrated density of states for a randomly perturbed lattice}

\author[R. Fukushima]{Ryoki Fukushima}

\address{%
Department of Mathematics\\ 
Kyoto University\\ 
Kyoto 606-8502\\ 
JAPAN\\
Current address:\\
Department of Mathematics\\
Tokyo Institute of Technology\\ 
Tokyo 152-8551\\
 JAPAN
}

\email{ryoki@math.titech.ac.jp}

\thanks{The first author was partially supported by JSPS Fellowships for Young Scientists.\\
The second author was partially supported by  KAKENHI (21540175)}
\author[N. Ueki]{Naomasa Ueki}
\address{Graduate School of Human and Environmental Studies\\
Kyoto University\\
Kyoto 606-8501\\ 
JAPAN} 
\email{ueki@math.h.kyoto-u.ac.jp}




\begin{abstract}
The asymptotic behavior of the integrated density of states for a randomly perturbed lattice
at the infimum of the spectrum is investigated.
The leading term is determined when the decay of the single site potential is slow.
The leading term depends only on the classical effect from the scalar potential.
To the contrary, the quantum effect appears
when the decay of the single site potential is fast.
The corresponding leading term is estimated and the leading order is determined.
In the multidimensional cases, the leading order varies in different ways from the known results in the Poisson case.
The same problem is considered for the negative potential.
These estimates are applied to investigate the long time asymptotics of Wiener integrals associated 
with the random potentials.
\end{abstract}

\maketitle

\section{Introduction}
In this paper, we are concerned with the self-adjoint operator in the form of 
\begin{equation}
 H_{\xi}=-h\Delta + \sum_{q \in \Z^d} u(\,\cdot-q-\xi_q)\label{operator}
\end{equation}
defined on the $L^2$-space on $\R ^d\setminus \bigcup _{q\in \Z ^d}(q+\xi _q+K)$
with the Dirichlet boundary condition, 
where $h$ is a positive constant and $K$ is a compact set in $\R ^d$ allowed to be empty. 
Our assumptions on the potential term are the following: 
(i) $\xi =(\xi_q)_{q \in \Z^d}$ is a collection of independent and identically
distributed $\R^d$-valued random variables with 
\begin{equation}
 \P(\xi_q \in dx) = \exp (-|x|^{\theta} )dx/Z(d,\theta) \label{single-prob} 
\end{equation}
for some $\theta > 0$ and the normalizing constant $Z(d,\theta)$; 
(ii) $u$ is a nonnegative function belonging to the Kato class $K_d$ (cf.\ \cite{CL90} p-53) and satisfying
\begin{equation}
u(x)=C_0|x|^{-\alpha}(1+o(1))\label{single-decay}
\end{equation} 
as $|x|\to \infty$ for some $\alpha >d$ and $C_0>0$. 

Although we assume the equality in \eqref{single-prob}, it will be easily seen from the proofs 
that only the asymptotic relation 
$$
	\P (\xi _q \in x+[0,1]^d) \asymp \exp (-|x|^{\theta})
$$
is essential for our theory, where $f(x) \asymp g(x)$ means 
$0<\varliminf _{|x| \to \infty}f(x)/g(x)\le 
\varlimsup _{|x| \to \infty}f(x)/g(x)<\infty$. 
In particular, we may replace $|x|^{\theta}$ by $(1+|x|)^{\theta}$ in \eqref{single-prob}.
Then the point process $\{ q+\xi _q\} _{q\in \Z ^d}$ converges weakly to
the complete lattice $\Z ^d$ as $\theta \to \infty$. Moreover, it is shown in Appendix A of~\cite{Fuk09a}
that this point process converges weakly to the Poisson point process with the intensity 1 as 
$\theta \downarrow 0$.
Since the Poisson point process is usually regarded as a completely disordered configuration, 
our model gives an interpolation between complete lattice and completely disordered media.

We will consider the integrated density of states $N(\lambda )$ ($\lambda \in \R$) of $H_{\xi}$ 
defined by the thermodynamic limit
\begin{equation}
\frac{1}{|\Lambda_R|} N_{\xi , \, \Lambda_R}(\lambda )\longrightarrow N(\lambda )\qquad 
\text{ as }R\to \infty. \label{m:ids}
\end{equation}
In \eqref{m:ids} we denote by $\Lambda _R$ a box $(-R/2,R/2)^d$ and by $N_{\xi ,\, \Lambda_R}(\lambda )$ 
the number of 
eigenvalues not exceeding $\lambda$ of the self-adjoint operator $H_{\xi , R}^D$ defined by 
restricting $H_{\xi}$ to $\Lambda _R\setminus \bigcup _{q\in \Z ^d}(q+\xi _q+K)$ 
with the Dirichlet boundary condition. 
We here note that the potential term in \eqref{operator} belongs to the local Kato class $K_{d,loc}$ 
(cf.\ \cite{CL90} p-53) as we will show in Section \ref{Appendix} below. 
It is then well known that the above limit exists for almost every $\xi$ and defines
a deterministic increasing function $N(\lambda)$
(cf.\ \cite{CL90}, \cite{KM82a}). 

The following are first two main results in this paper. 
\begin{thm}\label{pastur}
If $d<\alpha \le d+2$ and 
\begin{equation}
\essinf _{|x|\le R}u(x) \text{ is positive for any }R\ge 1, \label{positive}
\end{equation} 
then we have 
\begin{equation} 
 \log N(\lambda ) \asymp -\lambda ^{-\kappa}, \label{pastur-1}
\end{equation}
where $\kappa =(d+\theta )/(\alpha -d)$, and $f(\lambda ) \asymp g(\lambda )$ means $0<\varliminf _{\lambda \downarrow 0}f(\lambda )/g(\lambda )\le \varlimsup _{\lambda \downarrow 0}f(\lambda )/g(\lambda )<\infty$.
Moreover if $\alpha <d+2$, then we have
\begin{equation} 
 \lim_{\lambda \downarrow 0} \lambda ^{\kappa}\log N(\lambda )
=\frac{-\kappa ^{\kappa}}{(\kappa +1)^{\kappa +1}}\Bigg\{ \int_{\R ^d}dq \inf_{y\in \R ^d}\Big( \frac{C_0}{|q+y|^{\alpha}}+|y|^{\theta}\Big) \Bigg\} ^{\kappa +1}, \label{pastur-2}
\end{equation}
where the right hand side is finite by the assumption $\alpha >d$.
\end{thm}

\begin{thm}\label{lifshitz}
If $d=1$ and $\alpha>3$, then we have
 \begin{equation}
   \lim_{\lambda \downarrow 0}\lambda ^{(1+\theta )/2}\log N(\lambda) =
   -\frac{\pi ^{1+\theta}h^{(1+\theta )/2}}{(1+\theta )2^{\theta}}. \label{lifshitz-1}
 \end{equation}
If $d=2$ and $\alpha>4$, then we have
 \begin{equation}
   \log N(\lambda) \asymp
   -\lambda^{-1-\theta /2}\left( \log \frac{1}{\lambda} \right)^
    {-\theta /2}. \label{lifshitz-2}
 \end{equation}
If $d \ge 3$ and $\alpha > d+2$, then we have 
 \begin{equation}
   \log N(\lambda) \asymp -\lambda^{-(d+\mu \theta )/2}, \label{lifshitz-3}
 \end{equation}
where $\mu =2(\alpha -2)/(d(\alpha -d))$. 
\end{thm}

These results are generalizations of Corollary 3.1 in \cite{Fuk09a} to the case
that $\supp (u)$ is not compact (cf.\ Theorem \ref{cpt} below).
The results in Theorem \ref{pastur} are independent of the constant $h$.
In fact these asymptotics coincide with those of the corresponding classical integrated
density of states defined by 
$$N_c(\lambda )=\E [|\{ (x,p)\in \Lambda _R\times \R ^d : H_{\xi ,c}(x,p)\le \lambda \} |](2\pi \sqrt{h}R)^{-d}$$
for any $R\in \N$, where $| \cdot |$ is the $2d$-dimensional Lebesgue measure and
$$H_{\xi ,c}(x,p)=\sum_{j=1}^dp_j^2+V_{\xi}(x)$$
is the classical Hamiltonian (cf.\ \cite{LMW03}).
Therefore we may say that only the classical effect from the scalar potential determines 
the leading term for $\alpha <d+2$ and the leading order for $\alpha \le d+2$. 
To the contrary, the right hand side of \eqref{lifshitz-1} depends on $h$ and the right hand sides of \eqref{lifshitz-2} and \eqref{lifshitz-3} are strictly less than that of \eqref{pastur-1}.
Therefore we may say that the quantum effect appears in Theorem \ref{lifshitz}.
We here note that the right hand side of \eqref{pastur-1} gives an upper bound and the asymptotics of the classical counterpart not only for $\alpha \le d+2$ but also for $\alpha >d+2$ (see Proposition \ref{pastur-upper} below). 
For the critical case $\alpha =d+2$, the quantum effect appears at least in some cases.
We shall elaborate on this in Section \ref{sec-critical} below. 

In our model, the single site potentials are randomly displaced from the lattice. 
As is mentioned in \cite{Fuk09a}, such a model describes the Frenkel disorder in solid state physics 
and is called the random displacement model in the theory of random Schr\"odinger operator. 
Despite of the appropriateness of this model in physics, there are only a few mathematical studies and 
in particular the displacements have been assumed to be bounded in almost all works. 
For that case, Kirsch and Martinelli \cite{KM82b} discussed the existence of band gaps 
and Klopp \cite{Klo93} proved spectral localization in a semi-classical limit. 
More recently, Baker, Loss and Stolz \cite{BLS07}, \cite{BLS08} studied 
which configuration minimizes the spectrum of \eqref{operator} and also showed that the corresponding integrated density of states increases rapidly at the minimum in a one-dimensional example.
On the other hand, our displacements are unbounded.
Then the infimum of the spectrum is easily shown to be 0 opposed to the bounded cases.
This is an essential condition for our method, by which we investigate the behavior of $N(\lambda )$ at $\lambda =0$.
All our results show that $N(\lambda )$ increases slowly. 

In a slightly broader class of models where the potentials are randomly located, 
the most studied model is the Poisson model, 
where the random points $(q+\xi_q)_{q \in \Z^d}$ are replaced by the sample points of 
the Poisson random measure (cf.\ \cite{CL90}, \cite{PF92}). 
In the limit of $\theta \downarrow 0$, the above results coincide with the corresponding results 
for the Poisson model obtained by Pastur \cite{Pas77}, Lifshitz \cite{Lif65},
Donsker and Varadhan \cite{DV75c}, Nakao \cite{Nak77}, and \^Okura \cite{Oku81}.
As in the Poisson model, the critical value is always $\alpha =d+2$ and, in the one-dimensional case,
the leading order increases continuously as $\alpha$ increases to $d+2$ and does not depend on $\alpha \ge d+2$.
However in contrast to the Poisson case, the leading order jumps at $\alpha =d+2$ for $d=2$, and it 
depends on $\alpha \ge d+2$ for $d\ge 3$. 
These phenomena are due to the fact that the effect from states which have many tiny holes including $\{ q+\xi _q\} _q$ 
in their supports appears in the leading term of the asymptotics, as observed in \cite{Fuk09a}.
This is a characteristic difference with the Poisson case.
On the other hand, the decay rates of $N(\lambda )$ explode in the limit $\theta \to \infty$. 
This reflects the fact that the infimum of the spectrum is positive in the case of a finitely perturbed lattice including the case of the unperturbed lattice.

On the subjects of this paper, we have more results for the alloy type model
$$H_{\omega}=-h\Delta +\sum_{q\in \Z ^d}\omega _qu(x-q)$$
and the same critical value $\alpha =d+2$ is obtained,
where $\omega =(\omega _q)_{q\in \Z ^d}$ is a collection of independent and identically
distributed nonnegative real valued random variables.
As for the results, further developments and the relation with other models, refer to a recent survey by Kirsch and Metzger \cite{KM07}.

Our proof of Theorem \ref{pastur} is an extension of that of the corresponding result for the Poisson case (cf.\ \cite{Pas77}, \cite{PF92}).
For the proof of the multidimensional results in Theorem \ref{lifshitz}, we use a method based on a functional analytic approach (cf.\ \cite{CL90}, \cite{KM82a}).
This is different from the method in \cite{Fuk09a}, where a coarse graining method following Sznitman \cite{Szn98} is applied.
The method employed here can also be used to give a simpler proof of the results in the compact case in \cite{Fuk09a}.
We will present it in Section \ref{sec-lifshitz} below.
For the 1-dimensional result, we use a simple effective estimate of the first eigenvalue in \cite{Szn98}. 

As an application, we study the survival probability of the Brownian motion in a random environment. 
This was the main motivation in \cite{Fuk09a}.
We recall the connection between this and the integrated density of states, and extend the theory to the present settings. 
For the results, see Theorem \ref{F-K} below.
In the proof, we take the \emph{hard obstacles} $K$
appropriately so that the local singularity of the potential $u$ does not bring difficulty.
This is our only motivation to introduce the hard obstacles, and the hard obstacles do not affect the results. 

We also consider the operator
\begin{equation}
 H_{\xi}^-=-h\Delta - \sum_{q \in \Z^d} u(\,\cdot-q-\xi_q)\label{negative-operator}
\end{equation}
obtained by replacing the potential $u$ in $H_{\xi}$ by $-u$. 
For this operator, we assume $K=\emptyset$ since we are interested only in the effect of the negative potential.
The spectrum of this operator extends to $-\infty$.
For the asymptotic distribution, we show the following:

\begin{thm}\label{negative}
Suppose $K=\emptyset$, $\sup u=u(0)<\infty$ and $u(x)$ is lower semicontinuous at $x=0$. 
Then the integrated density of states $N^-(\lambda )$ of $H_{\xi}^-$ satisfies 
\begin{equation} 
 \lim_{\lambda \downarrow -\infty} \frac{\log N^-(\lambda )}{(-\lambda )^{1+\theta /d}}=\frac{-C_1}{u(0)^{1+\theta /d}},
 \label{negative-eq}
\end{equation}
where $C_1=d^{1+\theta /d}/\{ (d+\theta )|S^{d-1}|^{\theta/d}\}$ and $|S^{d-1}|$ is the volume of the $(d-1)$-dimensional surface $S^{d-1}$. 
\end{thm}

For the Poisson model, Pastur \cite{Pas77} showed that the corresponding
integrated density of states $N_{\rm Poi}^-(\lambda )$ satisfies
$$ \lim_{\lambda \downarrow -\infty} \frac{\log N_{\rm Poi}^-(\lambda )}{(-\lambda )\log (-\lambda )}=\frac{-1}{u(0)}.$$
The power of $\lambda$ in \eqref{negative-eq} tends to that of the Poisson model as $\theta \downarrow 0$. 
However, the logarithmic term is not recovered. Therefore, we cannot \emph{interchange the limits} 
$\lambda \downarrow -\infty$ and $\theta \downarrow 0$ in this case. 
Both for the Poisson and our cases, only the classical effect from the scalar potential determines the leading terms.
The lower semicontinuity of $u$ at 0 is a sufficient condition for the classical behavior: by this condition, the tunneling effect is suppressed.
For this subject, refer to Klopp and Pastur \cite{KP99}. 

Let us briefly explain the organization of this paper. 
We prove Theorems \ref{pastur}, \ref{lifshitz}, and \ref{negative} in Sections \ref{sec-pastur}, 
\ref{sec-lifshitz}, and \ref{sec-negative}, respectively. 
In Section \ref{sec-lifshitz} we also give a simple proof of the corresponding results for the case that $\supp (u)$ is compact.
In Section \ref{sec-critical}, we discuss the critical case $\alpha =d+2$.
In Section \ref{sec-FK} we study the asymptotic behaviors of certain Wiener integrals. 

\section{Proof of Theorem \ref{pastur}} \label{sec-pastur}
\subsection{Upper estimate}
To derive the asymptotics of the integrated density of states, one of the standard ways is to estimate 
its Laplace transform and use the Tauberian theorem (cf.~\cite{Fuk73,Nak77}).
We here say the Tauberian theorem by the theorem deducing the asymptotics from that of the Laplace-Stieltjes transform. 
Let $\NT (t)$ be the Laplace-Stieltjes transform of the integrated density of states $N(\lambda )$:
$$\NT (t)=\int_0^{\infty}e^{-t\lambda}dN(\lambda ).$$
Then, in view of the exponential Tauberian theorem due to Kasahara~\cite{Kas78}, 
the proof of the upper bound is reduced to the following:

\begin{prop}\label{pastur-upper}
If $K=\emptyset$ and \eqref{positive} is satisfied, then we have
 \begin{equation}
  \mathop{\varlimsup}_{t\uparrow \infty} \frac{\log \NT (t)}{t^{(d+\theta )/(\alpha +\theta )}}\le -\int_{\R ^d}dq \inf_{y\in \R ^d}\Big( \frac{C_0}{|q+y|^{\alpha}}+|y|^{\theta}\Big) 
 \end{equation}
for any $\alpha >d$.
\end{prop}

\noindent {\it Proof}.
We use the bound
 \begin{equation}
  \NT (t)\le \NT _1(t)(4\pi t h)^{-d/2},\label{upper-1}
 \end{equation}
where
$$\NT _1(t) = \int_{\Lambda _1}dx\E\Bigg[ \exp \Big( -t\sum_{q \in \Z^d} u(\,x-q-\xi_q)\Big) \Bigg] .$$
This is a simple modification of the bound in Theorem (9.6) in \cite{PF92}
for $\Z ^d$-stationary random fields. 
By replacing the summation by integration, we have
$$\log \NT _1(t)\le \int _{\R ^d}dq\log \E\Big[ \exp \Big( -t \inf_{x\in \Lambda _2}u(x-q-\xi _0)\Big) \Big] .$$ 
We pick an arbitrary $L>0$ and restrict the integration to $|q|\le Lt^{\eta}$. 
The assumption \eqref{single-decay} tells us that for any $\e _1>0$, there exists $R_1$ such that
$u(x)\ge C_0(1-\e _1)|x|^{-\alpha} \text{ whenever }|x|_{\infty}\ge R_1$, 
where $|x|_{\infty}=\max_{1\le i\le d}|x_i|$.
Thus the right hand side is dominated by
 \begin{equation*}
  \begin{split}
   & \int _{|q|\le {Lt^{\eta}}}dq\log \Bigg\{ \int _{|q+y|_{\infty}\ge R_1+1}
\frac{dy}
{Z(d,\theta )}\exp \Big( -t \inf_{x\in \Lambda _2}\frac{C_0(1-\e _1)}{|x-q-y|^{\alpha}}-|y|^{\theta}\Big) \\
& +\exp \Big( -t\inf_{\Lambda _{2R_1+4}}u\Big) \Bigg\} .
  \end{split} 
 \end{equation*}
Thanks to the assumption \eqref{positive}, the second term makes only negligible contribution to the asymptotics. 
By changing the variables $(q,y)$ to $(t^{-\eta}q, t^{-\eta}y)$ with $\eta =1/(\alpha +\theta )$, 
we see that this equals
$$t^{d\eta}\int _{|q|\le L}dq\log \Big\{ \NT_2(t, q)+\exp \Big( -t\inf_{\Lambda _{2R_1+4}}u\Big) \Big\} ,$$
where 
 \begin{equation*}
  \begin{split}
& \NT _2(t, q) \\
& = t^{d\eta }\int _{|q+y|_{\infty}\ge (R_1+1)t^{-\eta }}
\frac{dy}
{Z(d,\theta )}\exp \Big( -t^{\theta \eta } \inf_{x\in \Lambda _{2t^{-\eta }}}\frac{C_0(1-\e _1)}
{|x-q-y|^{\alpha}}-t^{\theta \eta }|y|^{\theta}\Big).
  \end{split} 
 \end{equation*}
We take $L$ as an arbitrary constant independent of $t$. 
Then, taking $\e _2, \e _3>0$ sufficiently small, 
we can dominate $\NT _2(t, q)$ by
$\exp ( -t^{\theta \eta } \NT _3(q))\e _2^{-d/\theta}$ for large enough $t$, where
$$\NT _3(q) = \inf \Big\{ \frac{C_0(1-\e _1)}{|x-q-y|^{\alpha}}+(1-\e _2)|y|^{\theta} 
: x\in \Lambda _{\e _3}, y\in \R ^d\Big\} .$$
Therefore we obtain
$$\mathop{\varlimsup}_{t\uparrow \infty} \frac{\log \NT (t)}{t^{(d+\theta )\eta }}\le -\int_{|q|\le L}\NT _3(q)dq.$$
Since $\e _1, \e _2, \e _3$ and $L$ are arbitrary,
this completes the proof.
\qed

\subsection{Lower estimate}
To prove the lower estimate, we have only to show the following:

\begin{prop}\label{pastur-lower}
If $\alpha <d+2$, then we have
 \begin{equation}
  \mathop{\varliminf}_{t\uparrow \infty} \frac{\log \NT (t)}{t^{(d+\theta )/(\alpha +\theta )}}\ge -\int_{\R ^d}dq \inf_{y\in \R ^d}\Big( \frac{C_0}{|q+y|^{\alpha}}+|y|^{\theta}\Big).
 \end{equation}
Moreover, this bound remains valid for $\alpha = d+2$ with a smaller constant in the right hand side. 
\end{prop}

The case $\alpha = d+2$ will be discussed in more detail in Section \ref{sec-critical} below.

\noindent {\it Proof of Proposition \ref{pastur-lower}}.
We use the bound
 \begin{equation}
  \NT (t)\ge R^{-d}\exp (-th\| \nabla \psi _R\| _2^2)\NT _1(t)\label{lower-1}
 \end{equation}
which holds for any $R\in \N$ and $\psi _R\in C_0^{\infty} (\Lambda _R)$ such that 
$\| \psi _R\| _2=1$, where $\| \cdot \| _2$ is the $L^2$-norm, and
$$\NT _1(t) = \E\Bigg[\exp \Big( -t\sum_{q \in \Z^d} \int dx \psi _R(x)^2u(\,x-q-\xi_q)\Big) : 
\bigcup _{q\in \Z ^d}(q+\xi _q+K)\cap \Lambda _R=\emptyset \Bigg] .$$
This can be proven by the same method as for the corresponding bound in Theorem (9.6) in \cite{PF92}
for $\R ^d$-stationary random fields.
By replacing the summation by integration, we have
$$\log \NT _1(t)\ge \int _{\R ^d}\NT _2(t,q)dq,$$
where
 \begin{equation*}
  \begin{split}
\NT_2(t,q) 
= \log \E\Big[ & \exp \Big( -t \int dx \psi _R(x)^2\sup_{z\in \Lambda _1}u(x-q-z-\xi _0)\Big) \\
& : (q+\xi _0+K)\cap \Lambda _R=\emptyset \Big] .
  \end{split} 
 \end{equation*}
For any $\e _1>0$, there exists $R_1$ such that $K\subset B(R_1)$ and $u(x)\le C_0(1+\e _1)|x|^{-\alpha}$ 
for any $|x|\ge R_1$ 
by the assumption \eqref{single-decay}.
To use this bound in the above right hand side, we need 
$\inf \{ |x-q-z-\xi_0| : x\in \Lambda _R, z\in \Lambda _1\} \ge R_1$. 
However we shall deal with a simpler sufficient condition $|\xi_0| \le |q|/2$ and $|q|\ge 2(R_1+\sqrt{d}R)$ instead.
Now fix $\beta >0$ and take $t$ large enough so that $t^{\beta} >2(R_1+\sqrt{d}R)$. 
Then we obtain
\begin{equation}
\int _{|q|\ge t^{\beta}}\NT _2(t,q)dq
\ge \int _{|q|\ge t^{\beta}}dq\Big( -\frac{tC_0(1+\e _1)2^{\alpha}}{(|q|-2\sqrt{d}R)^{\alpha}}+\log \P(|\xi _0|\le |q|/2)\Big) .\label{out}
\end{equation}
By a simple estimate using 
$\log (1-X)\ge -2X$ for $0\le X\le 1/2$, we can bound the right hand side
from below by
$-c_1t^{1-\beta (\alpha -d)}-c_2\exp (-c_3t^{\beta \theta})$.
The other part is estimated as
 \begin{equation}
  \begin{split}
   & \int _{|q|\le t^{\beta}}\NT _2(t,q)dq \\
   \ge & \mathop{\int}_{|q|\le t^{\beta}}dq\log \mathop{\int}_{|q+y|\ge R_1+\sqrt{d}R}\frac{dy}{Z(d,\theta )} \\
   & \times \exp \Big( -\frac{tC_0(1+\e _1)}{\inf \{ |x-q-z-y|^{\alpha} : x\in \Lambda _R, z\in \Lambda _1\}}-|y|^{\theta}\Big) .
  \end{split} \label{lower-2}
 \end{equation}
By changing the variables, we find that the right hand side equals
$$t^{d\eta }\mathop{\int}_{|q|\le t^{\beta -\eta }}dq\log \mathop{\int}_{|q+y|
\ge (R_1+\sqrt{d}R)t^{-\eta }}\frac{dyt^{d\eta }}{Z(d,\theta )}\exp (-t^{\theta \eta }\NT _3(y,q)),$$
where $\eta =1/(\alpha +\theta )$ and
\begin{equation}
\NT _3(y,q) = \frac{C_0(1+\e _1)}{\inf \{ |x-q-z-y|^{\alpha} : x\in \Lambda _{Rt^{-\eta }},z\in \Lambda _{t^{-\eta }}\}}+|y|^{\theta}. \label{n3t}
\end{equation}
Let us take $\gamma >0$ and restrict the integration with respect to $y$ to the ball $B(y_0,t^{-\gamma})$ 
with center $y_0$ and radius $t^{-\gamma}$.
Then we can bound the integrand with respect to $q$ from below by
\begin{equation}
\log \frac{|B(0,1)| t^{d(\eta -\gamma )}}{Z(d,\theta )}-t^{\theta \eta }\NT _4(q,t),\label{lower-3}
\end{equation}
where  
\begin{equation}
  \begin{split}
   \NT _4(q,t) = \inf \Big\{ & \mathop{\sup}_{y\in B(y_0,t^{-\gamma})}\NT _3(y,q) \\
& : y_0\in \R ^d, d(B(y_0,t^{-\gamma}),-q)\ge (R_1+\sqrt{d}R)t^{-\eta }\Big\} .
  \end{split} \label{n4t}
\end{equation}
We now specify $R$ as the integer part of $\e _2t^\eta $, where $\e _2$ is an arbitrarily fixed
positive number. We take $\psi _R$ as the nonnegative and normalized ground state of the Dirichlet Laplacian on the cube 
$\Lambda _R$ and take $\beta$ between $\eta $ and $\eta (1+\theta /d)$.
Then, for $\alpha <d+2$, we obtain
\begin{equation}
\mathop{\varliminf}_{t\uparrow \infty} \frac{\log \NT (t)}{t^{(d+\theta )\eta }}
\ge -\mathop{\varlimsup}_{t\uparrow \infty} \int_{|q|\le t^{\beta -\eta }}dq \NT _4(q,t),
\label{lower-4}
\end{equation}
since $th\| \nabla \psi_R \|_2 \asymp tR^{-2}$ and \eqref{out} is negligible compared with $t^{(d+\theta )\eta }$. 
When $|q|\le t^{\beta -\eta }$, we can dominate $1/t$ by a power of $q$.
Thus, for large $|q|$, by taking $y_0$ as $0$, we can dominate $\NT _4(q,t)$ by
$|q|^{-\alpha}+|q|^{-\gamma \theta /(\beta -\eta )}$. 
This is integrable if we take $\gamma$ large enough so that $\gamma \theta /(\beta -\eta )>d$.
Thus, by the Lebesgue convergence theorem, we have
 \begin{equation*}
  \begin{split}
& \lim_{t\uparrow \infty} \int_{|q|\le t^{\beta -\eta }}dq \NT _4(q,t)\\
& =\int_{\R ^d}dq\inf \Big\{ \frac{C_0(1+\e _1)}
{\ds \inf_{x\in \Lambda _{\e _2}}|x-q-y|^{\alpha}}+|y|^{\theta} : y\in \R ^d, d(y,q)\ge \e _2\sqrt{d}\Big\} .
  \end{split} 
 \end{equation*}
Since $\e _1$ and $\e _2$ are arbitrary,
this completes the proof of the former part of Proposition \ref{pastur-lower}. 
For the case $\alpha = d+2$, we take $\e_2 = 1$. Then we have $th\| \nabla \psi_R \|_2 \asymp t^{(d+\theta )\eta }$ 
and the latter part of Proposition \ref{pastur-lower} follows from the same argument as above. 
\qed

\section{Proof of Theorem \ref{lifshitz}~and the compact case} \label{sec-lifshitz}
In this section, we use some additional notations to simplify the  presentation. 
For any self-adjoint operator $A$, let $\lambda _1(A)$ be the infimum of its spectrum and, 
for any locally integrable function $V$ and $R>0$, let $(-h\Delta +V)_R^D$ and $(-h\Delta +V)_R^N$ be the self-adjoint operators 
$-h\Delta +V$ on the $L^2$-space on the cube $\Lambda _R$ with the Dirichlet and the Neumann boundary conditions, respectively.
\subsection{Proof of Theorem \ref{lifshitz} (I): One-dimensional case}
To obtain the upper estimate, we have only to show the following:

\begin{prop}\label{1D-upper}
If $d=1$, $K=\emptyset$, $\supp (u)$ is compact, 
\begin{equation}
\liminf_{x\downarrow 0}\int_0^xu(y)dy/x>0 \text{, and }\liminf_{x\downarrow 0}\int_{-x}^0u(y)dy/x>0, \label{1D-upper-ass}
\end{equation} 
then we have
 \begin{equation}
   \mathop{\varlimsup}_{t\uparrow \infty} \frac{\log \NT (t)}{t^{(1+\theta )/(3+\theta )}} \le
    -\frac{3+\theta}{1+\theta}\Big( \frac{h\pi ^2}4 \Big) ^{(1+\theta )/(3+\theta )}.
 \label{1D-upper1}
 \end{equation}
\end{prop}

\noindent {\it Proof}.
We assume $h=1$ for simplicity.
In the well known expression 
$$\NT (t)=\int_{\Lambda _1}\E[\exp (-tH_{\xi})(x,x)]dx,$$
we apply the Feynman-Kac formula and an estimate for the exit time of the Brownian motion (cf.\ \cite{IM74}) to obtain
$$\NT (t)\le \int_{\Lambda _1}\E[\exp (-tH^D_{\xi ,t})(x,x)]dx+c_1e^{-c_2t},$$
where $\exp (-tH_{\xi})(x,y)$ and $\exp (-tH^D_{\xi ,t})(x,y)$, $t>0$,
$x, y\in \R$, are the integral kernels of the heat semigroups generated by $H_{\xi}$ and $H^D_{\xi ,t}$, respectively.
By the eigenfunction expansion of the integral kernel, we have
$$\NT (t)\le c_3t\NT _1(t)+c_4e^{-c_5t},$$
where $\NT _1(t)=\E[\exp (-t\lambda _1(H^D_{\xi ,t}))]$.
Thus we have only to prove \eqref{1D-upper1} with $\NT (t)$ replaced by $\NT _1(t)$.
Now we use Theorem 3.1 in the page 123 in \cite{Szn98}, which states 
$$\lambda _1(H_{\xi ,t}^D)\ge \pi ^2/(\sup _k|I_k|+c_6)^2$$
for large enough $t$ under the assumption \eqref{1D-upper-ass}, 
where $\{ I_k\} _k$ are the random open intervals such that $\sum _kI_k=\Lambda _t-\{ q+\xi _q : q\in \Z \}$ and 
$|I_k|$ is the length of $I_k$.
If $\sup _k|I_k|\ge s$ for some $0\le s \le t$, then there exists $p\in \Z \cap \Lambda _t$ such that $\{ q+\xi _q : q\in \Z\} \cap [p,p+s-2] =\emptyset$.
The probability of this event is estimated as
 \begin{equation*}
  \begin{split}
   & \P(\sup _k|I_k|\ge s)\le \sum_{p\in \Z \cap \Lambda _t}\prod_{q\in \Z \cap [p,p+s-2]}\P(q+\xi _q\not\in [p,p+s-2]) \\
   & \le t\prod_{q\in \Z \cap [p,p+s-2]}\exp (-(1-\varepsilon )d(q,[p,p+s-2]^c)^{\theta})/\varepsilon ^{1/\theta} \\
   & \le t\exp \Big( -(1-\varepsilon )\int _0^{s-3} d(q,[0,s-3]^c)^{\theta}dq+\frac{s}{\theta}\log \frac1{\varepsilon}\Big) \\
   & \le t\exp \Big( -\frac{2(1-\varepsilon )}{\theta +1}\Big( \frac{s-3}2\Big) ^{\theta +1}+\frac{s}{\theta}\log \frac1{\varepsilon}\Big)
  \end{split} 
 \end{equation*}
if $s\ge 3$, where $0<\varepsilon <1$ is arbitrary.
Therefore we have
 \begin{equation*}
   \NT _1(t) 
   \le c_7t^2\exp \Big( -\inf _{R>3} \Big( t\frac{\pi^2}{(R+c_6)^2}
   +\frac{(1-\varepsilon )}{2^{\theta}(\theta +1)} ( R-3 )^{\theta +1}
   -\frac{R}{\theta} \log \frac1{\varepsilon}\Big) \Big)  
   +c_8e^{-c_9t}
 \end{equation*}
for large $t$.
Now it is easy to see that the infimum in the right hand side is attained 
by $R \sim 2(\pi^2 t/4)^{1/(3+\theta)}$ and we obtain \eqref{1D-upper1}.
\qed
\begin{rem}
We put the additional assumption \eqref{1D-upper-ass}
only to use Theorem 3.1 in the page 123 in \cite{Szn98}. These assumptions are not restrictive at all 
since we can always find a $z \in \R$ such that $u(\,\cdot\,+z)$ satisfies them 
by the fundamental theorem of calculus 
and such a finite translation of $u$ does not affect the above argument. 
\end{rem}

\begin{prop}\label{1D-lower}
If $d=1$ and $\alpha >3$, then we have
 \begin{equation}
   \mathop{\varliminf}_{t\uparrow \infty} \frac{\log \NT (t)}{t^{(1+\theta )/(3+\theta )}} \ge
    -\frac{3+\theta}{1+\theta}\Big( \frac{h\pi ^2}4 \Big) ^{(1+\theta )/(3+\theta )}.
 \label{1D-lower1}
 \end{equation}
\end{prop}

\noindent {\it Proof}.
This is proven by modifying our proof of Proposition \ref{pastur-lower}.
We take $\psi _R$ as the nonnegative and normalized ground state of $(-\Delta )^D_R$.
In \eqref{lower-2}, we restrict the integral with respect to $y$ to $|q+y|\ge R_1+(R+1)/2$.
In \eqref{lower-3}, we take $\eta =1/(3+\theta )$ and $R$ as the integer part of ${\mathcal R}t^{\eta}$ for a positive number ${\mathcal R}>0$.
Then since $t\| \nabla \psi_R \|_2^2 \sim t^{(1+\theta )\eta}(\pi /{\mathcal R})^2$ is not negligible,
\eqref{lower-4} is modified as
$$\mathop{\varliminf}_{t\uparrow \infty} \frac{\log \NT (t)}{t^{(1+\theta )\eta }}
\ge -h\Big( \frac{\pi}{\mathcal R}\Big) ^2-\mathop{\varlimsup}_{t\uparrow \infty} \int_{|q|\le t^{\beta -\eta }}dq \NT _4(q,t),$$
where $\NT _4(q,t)$ is defined by replacing $\NT _3(y,q)$ and $R_1+\sqrt{d}$ by
$$\frac{C_0(1+\e _1)}{t^{(\alpha -3)\eta}\inf \{ |x-q-z-y|^{\alpha} : x\in \Lambda _{Rt^{-\eta }},z\in \Lambda _{t^{-\eta }}\}}+|y|^{\theta}$$
and $R_1+(R+1)/2$, respectively, in \eqref{n4t}.
Since
$$\mathop{\varlimsup}_{t\uparrow \infty}\NT _4(q,t)\le \inf_{y\not\in \Lambda _{\mathcal R}(-q)}|y|^{\theta}=d(q,\Lambda _{\mathcal R}^c)^{\theta},$$
we obtain
$$\mathop{\varliminf}_{t\uparrow \infty} \frac{\log \NT (t)}{t^{(1+\theta )\eta }}
\ge -h\Big( \frac{\pi}{\mathcal R}\Big) ^2-\frac{{\mathcal R}^{\theta +1}}{2^{\theta}(\theta +1)},$$
by the Lebesgue convergence theorem.
By taking the supremum over ${\mathcal R} >0 $, we obtain the result.
\qed

\subsection{Proof of Theorem \ref{lifshitz} (II) : Upper estimate for the multidimensional case}
In the two-dimensional case, we can simply use Corollary 3.1 in \cite{Fuk09a} to get the upper bound. 
Indeed, the integrated density of states increases if we truncate the tail of $u$ 
and hence the bound for the compactly supported potentials yields
 \begin{equation}
   N(\lambda ) \le
     c_1\exp ( -c_2\lambda ^{-1-\theta /2}(\log (1/\lambda ))^{-\theta /2})\label{two} ,
 \end{equation}
for $0\le \lambda \le c_3$, where $c_1$, $c_2$ and $c_3$ are positive constants depending on $h$ and $C_0$.
We give another proof for Corollary 3.1 in \cite{Fuk09a} in Subsection~\ref{compact} below.

In the rest of this subsection we assume $d\ge 3$.
Then our goal is the following: 

\begin{prop}\label{lifshitz-upper}
Let $\alpha \ge d+2$ and $K=\emptyset$. There exist finite positive function $k_1(h)$ and $k_2(h)$ of $h$ and a positive constant $c$ such that
 \begin{equation}
   N(\lambda) \le
    k_1(h)\exp (-c((h\wedge h^{(\alpha -d)/(\alpha -2)})/\lambda )^{(d+\mu \theta )/2}) \label{higher}
 \end{equation}
for $0\le \lambda \le k_2(h)$. 
\end{prop}

We first see that Proposition \ref{lifshitz-upper} follows from the following: 

\begin{prop}\label{ev-bdd}
For sufficiently small $\e _1, \e _2>0$,
there exist a positive constant $c$ independent of $(h, R)$, and positive constants $c'$ and $c''$ independent of $(c_0, h, R)$ such that $\# \{ q\in \Z ^d\cap \Lambda _R : |\xi _q|\ge \e _1R^{\mu}\} \le \e _2R^d$, $R^{\mu d}\ge c'h/c_0$ and $R^{\mu (\alpha -2-d)}\ge c''c_0/h$ imply
 \begin{equation}
  \lambda _1\Big( \Big( -h\Delta +\sum_{q\in \Z ^d\cap \Lambda _R}\frac{c_01_{B(q+\xi _q,R_0)^c}(x)}{|x-q-\xi _q|^{\alpha}}\Big) _R^N\Big)
  \ge c(h\wedge h^{(\alpha -d)/(\alpha -2)})/R^2,
 \end{equation}
where $c_0$ and $R_0$ are arbitrarily fixed positive constants and $1_D$ is the characteristic function of 
$D \subset \R ^d$.
\end{prop}

\noindent {\it Proof of Proposition \ref{lifshitz-upper}}.
It is well known that
$$N(\lambda )\le \frac{c_1}{(R\wedge \sqrt{h})^d}\P( \lambda _1(H_R^N) \le \lambda )$$
(cf.\ (10.10) in \cite{PF92}).
We can take $c_0$ and $R_0$ so that 
$$u(x)\ge c_01_{B(R_0)^c}(x)|x|^{-\alpha}.$$
Thus by Proposition \ref{ev-bdd}, there exists a constant $c_2$ such that
 \begin{equation*}
  \begin{split}
& N(c_2(h\wedge h^{(\alpha -d)/(\alpha -2)})/R^2)\\
& \le \frac{c_1}{(R\wedge \sqrt{h})^d}\P( \# \{ q\in \Z ^d\cap \Lambda _R : |\xi _q|\ge \e _1R^{\mu}\} \ge \e _2R^d).
  \end{split} 
 \end{equation*}
We here should take $c_0$ sufficiently small so that the conditions of Proposition \ref{ev-bdd} are satisfied if $\alpha =d+2$.
When the event in the right hand side occurs, we have
$$\sum_{q\in \Z ^d\cap \Lambda _R} |\xi _q|^{\theta}\ge \e _1^{\theta}\e _2R^{d+\mu \theta}.$$
Thus it is easy to show
$$N(c_2(h\wedge h^{(\alpha -d)/(\alpha -2)})/R^2)\le \frac{c_3}{(R\wedge \sqrt{h})^d}\exp (-c_4R^{d+\mu \theta}),$$
and \eqref{higher} follows immediately. 
\qed

We next proceed to the proof of Proposition \ref{ev-bdd}.
We start with the following:

\begin{lem}\label{unif-pt-1} 
 $\inf \{ \lambda _1((-\Delta +1_{B(b,1)})_R^N) : b\in \Lambda _R\} \ge cR^{-d}$. 
\end{lem}

This lemma follows immediately from the Proposition 2.3 of Taylor \cite{Tay76} using the scaling 
with the factor $R^{-1}$. 
That proposition is stated in terms of the scattering length.
We here give an elementary proof following a lemma in the page 378 in Rauch \cite{Rau75-LNM2} 
for the reader's convenience.
 
\begin{proof}
We rewrite as $\lambda _1((-\Delta +1_{B(b,1)})_R^N)=\lambda _1((-\Delta +1_{B(1)})_{R,b}^N)$, where, for any locally integrable function $V$ and $R>0$, $(-\Delta +V)_{R,b}^N$ is the self-adjoint operator $-\Delta +V$ on the $L^2$ space on the cube $\Lambda _R(b)=b+\Lambda _R$ with the the Neumann boundary condition, and $B(1)=B(0,1)$.
For any smooth function $\varphi$ on the closure of $\Lambda _R(b)$, we have
 \begin{equation*}
  \begin{split}
   & \int_{\Lambda _R(b)}\varphi ^2(x)dx\\
   & =\int_1^{R(b)}dr r^{d-1}\mathop{\int}_{\theta \in S^{d-1} : (r,\theta )\in \Lambda _R(b)}dS\Big( \varphi   (g(r),\theta )+\int_{g(r)}^r\partial _s\varphi (s,\theta )ds \Big) ^2\\
   & \quad +\int_{B(1)\cap \Lambda _R(b)}  
    \varphi ^2(x)dx,
  \end{split}
 \end{equation*}
where $(r,\theta )$ is the polar coordinate, $R(b)=\sup \{ |x| : x\in \Lambda _R(b)\}$, $dS$ is the volume element of the $(d-1)$-dimensional surface $S^{d-1}$ and $g(r)=\{ (r-1)/(R(b)-1)+1\} /2$.
By the Schwarz inequality and a simple estimate, we can show
$$\int_1^{R(b)}dr r^{d-1}\mathop{\int}_{\theta \in S^{d-1} : (r,\theta )\in \Lambda _R(b)}dS\Big( \int_{g(r)}^r\partial _s\varphi (s,\theta )ds \Big) ^2
\le cR(b)^d\int_{\Lambda _R(b)}|\nabla \varphi |^2(x) dx,$$
where $c$ is a constant depending only on $d$.
By changing the variable, we can also show
$$\int_1^{R(b)}dr r^{d-1}\mathop{\int}_{\theta \in S^{d-1} : (r,\theta )\in \Lambda _R(b)}dS\varphi   (g(r),\theta )^2
\le c'R(b)^d\int_{B(1)\cap \Lambda _R(b)}\varphi ^2(x) dx,$$
where $c'$ is also a constant depending only on $d$.
Since $\sup_{b\in \Lambda _R}R(b)\le \sqrt{d}R$, we can complete the proof.
\end{proof}

\begin{lem}\label{unif-pt-2} 
There exist positive constants $c$, $c'$, and $c''$ such that
 \begin{equation*}
  \begin{split}
& \inf \Big\{ \lambda _1\Big( \Big( -h\Delta +\sum_{j=1}^n\frac{c_01_{B(b_j,R_0)^c}(x)}{|x-b_j|^{\alpha}}\Big)_R^N\Big) : b_1,\ldots ,b_n\in \Lambda _R\Big\} \\
& \ge c(c_0n)^{(d-2)/(\alpha -2)}h^{(\alpha -d)/(\alpha -2)}/R^d
  \end{split}
 \end{equation*}
for $n \ge c'h/c_0$ and $R \ge c''(c_0n/h)^{1/(\alpha -2)}$. 
\end{lem}

\noindent {\it Proof}.
Since $\lambda _1(A+B)\ge \lambda _1(A)+\lambda _1(B)$ for any self-adjoint operators $A$ and $B$, 
the left hand side is bounded from below by
$$\inf \{ \lambda _1((-h\Delta +c_0n1_{B(b,R_0)^c}(x)|x-b|^{-\alpha})_R^N) : b\in \Lambda _R\}.$$
A change of the variable shows that this equals
$$hk^{-2}\inf \{ \lambda _1((-\Delta +c_0nk^{2-\alpha}h^{-1}1_{B(b,R_0/k)^c}(x)|x-b|^{-\alpha})_{R/k}^N) 
: b\in \Lambda _{R/k}\}$$ 
for any $k>0$. 
We can bound this from below by 
$$hk^{-2}\inf \{ \lambda _1((-\Delta +c_0nk^{2-\alpha}h^{-1}3^{-\alpha}1_{B(b',1)}(x))_{R/k}^N) 
: b'\in \Lambda _{R/k}\}$$ 
for $k\ge R_0$ and $R>4\sqrt{d}k$, and we can use Lemma \ref{unif-pt-1} to complete the proof by taking 
$k$ as $(c_0n3^{-\alpha}h^{-1})^{1/(\alpha -2)}$.
Indeed, for each $b\in \Lambda _{R/k}$, we set $b':=b-(1+R_0/k)b/|b|$ if $b$ is not the zero vector. 
If $b$ is the zero vector, we set $b'$ as an arbitrarily chosen vector with the norm $1+R_0/k$.
Since $R_0/k\le |x-b|\le 2+R_0/k$ on $B(b',1)$, we have 
$$1_{B(b,R_0/k)^c}(x)|x-b|^{-\alpha}\ge (2+R_0/k)^{-\alpha}1_{B(b',1)}(x).$$
We bound this from below by $3^{-\alpha}1_{B(b',1)}(x)$
by assuming $k\ge R_0$.
Moreover we claim  $b'\in \Lambda _{R/k}$ for all $b\in \Lambda _{R/k}$.
A sufficient condition for this is $R\ge 2\sqrt{d}(R_0+k)$,
since $b'$ for $b$ with $|b|\ge 1+R_0/k$ is a contraction of $b$ and $\sup \{ |b'|_{\infty} : |b|\le 1+R_0/k\}=\sqrt{d}(1+R_0/k)$.
\qed

\begin{lem}\label{unif-ef-bdd} 
Let $V$ be any locally integrable nonnegative function on $\R ^d$.
Then any eigenfunction $\phi$ of $(-h\Delta +V)_R^N$ satisfies
$$\| \phi \| _{\infty}\le c(1/R+\sqrt{\lambda /h})^{d/2}\| \phi \| _2,$$
where $c$ is a finite constant depending only on $d$, $\lambda$ is the corresponding eigenvalue, and $\| \cdot \| _{\infty}$ and $\| \cdot \| _2$ are $L^{\infty}$ and $L^2$ norms, respectively. 
\end{lem}

The proof of this lemma is same as that of (3.1.55) in \cite{Szn98}.
Now we prove Proposition \ref{ev-bdd}:

\noindent {\it Proof of Proposition \ref{ev-bdd}}.
We use the following classification:
$${\mathcal F}= \{ a\in \Lambda _R\cap R^{\mu}\Z ^d 
: \# (\Lambda _{R^{\mu}}(a)\cap \{ q+\xi _q : q\in \Z ^d\cap \Lambda _R\} )<R^{\mu d}/2\}$$
and
$${\mathcal N}= \{ a\in \Lambda _R\cap R^{\mu}\Z ^d 
: \# (\Lambda _{R^{\mu}}(a)\cap \{ q+\xi _q : q\in \Z ^d\cap \Lambda _R\} )\ge R^{\mu d}/2\} .$$
By Lemma \ref{unif-pt-2}, 
$$\lambda _1((-h\Delta +\sum_qc_01_{B(q+\xi _q,R_0)^c}(x)|x-q-\xi _q|^{-\alpha})_{R^{\mu},a}^N)\ge ch^{(\alpha -d)/(\alpha -2)}/R^2$$
for any $a\in {\mathcal N}$. 
Let us write $\varphi$ for the nonnegative and normalized ground state of the operator 
$(-h\Delta +\sum_qc_01_{B(q+\xi _q,R_0)^c}(x)|x-q-\xi _q|^{-\alpha}) _R^N$. 
Then, applying the Rayleigh--Ritz variational formula, we have
$$\lambda _1\Big( \Big( -h\Delta +\sum_q\frac{c_01_{B(q+\xi _q,R_0)^c}(x)}{|x-q-\xi _q|^{\alpha}}\Big) _R^N\Big) 
\ge \frac{ch^{(\alpha -d)/(\alpha -2)}}{R^2}\sum_{a\in {\mathcal N}}\int_{\Lambda _{R^{\mu}}(a)}\varphi ^2dx.$$
If we assume $\lambda _1((-h\Delta +\sum_qc_01_{B(q+\xi _q,R_0)^c}(x)|x-q-\xi _q|^{-\alpha})_{R^{\mu},a}^N)\le Mh/R^2$, then Lemma \ref{unif-ef-bdd} implies that the right hand side is bounded from below by
\begin{equation}
cR^{-2}h^{(\alpha -d)/(\alpha -2)}(1-c'M^{d/2}R^{(\mu -1)d}\# {\mathcal F}). \label{ev-bdd-pf1}
\end{equation}
Since $\# (\Lambda _{R^{\mu}}(a)\cap \{ q+\xi _q : q\in \Z ^d\cap \Lambda _R\} )
\ge \# \{ q\in \Lambda _{(1-2\e _1)R^{\mu}}(a)\cap \Z ^d : |\xi _q|\le \e _1R^{\mu}\}$, we have $\# \{ q\in \Lambda _{(1-2\e _1)R^{\mu}}(a)\cap \Z ^d : |\xi _q|\le \e _1R^{\mu}\} <R^{\mu d}/2$ and $\# \{ q\in \Lambda _{(1-2\e _1)R^{\mu}}(a)\cap \Z ^d : |\xi _q|\ge \e _1R^{\mu}\} >\{ (1-2\e _1)^d-1/2\} R^{\mu d}$ for $a\in {\mathcal F}$.
Thus, by the assumption of this proposition, we have
$\e _2R^d\ge (\# {\mathcal F})\{ (1-2\e _1)^d-1/2\} R^{\mu d}$ 
and $\# {\mathcal F}\le R^{d(1-\mu )}\e _2/\{ (1-2\e _1)^2-1/2\}$.
By substituting this to \eqref{ev-bdd-pf1}, we complete the proof.
\qed
\subsection{Proof of Theorem \ref{lifshitz} (III) : Lower estimate for the multidimensional case}
We shall work with $h=C_0=1$ for simplicity.

\begin{prop}\label{lifshitz-lower}
Suppose $d=2$ and $\alpha>4$ or $d \ge 3$ and $\alpha \ge d+2$. 
Then there exist positive constants $c_1$, $c_2$, and $c_3$ such that
 \begin{equation}
  \begin{split}
   N(\lambda) \ge
   \left\{
   \begin{array}{lr}
    \smash[t]{c_1\exp\Bigl(-c_2\lambda^{-1-\theta /2}\left( \log (1/\lambda ) \right)^
    {-\theta /2}\Bigr)} &(d=2),\\[8pt]
    c_1\exp (-c_2\lambda ^{-(d+\mu \theta )/2})  &(d \ge 3), 
   \end{array}\right.
  \end{split}
 \end{equation}
for $0\le \lambda \le c_3$. 
\end{prop}

\noindent {\it Proof}.
We consider the event
 \begin{equation}
  \begin{split}
   & \{ \text{For any }p\in R_1\Z ^d\cap \Lambda _{3R}\text{ and }q\in \Z ^d\cap \Lambda _{R_1}(p)\cap \Lambda_{2R}, q+\xi _q\in \Lambda _1(p) \} \\
   & \cap \{\text{For any }q\in \Z ^d \setminus \Lambda _{2R}, |\xi _q|\le |q|/4 \}
  \end{split} \label{lower-event}
 \end{equation}
where $R_1 = R^{\mu}$ for $d\ge 3$ and $R_1=R/\sqrt{\log R}$ for $d=2$. 
Then we have
\begin{equation}
  \begin{split}
   N(\lambda )\ge R^{-d}\P\Big( & \| \nabla \Phi _R\| _2^2+\Big( \Phi _R,\sum_{q\in \Z ^d}u(x-q-\xi _q)\Phi _R\Big) \le \lambda \\
   & \text{ and the event \eqref{lower-event} occurs}\Big),
  \end{split} \label{R-R}
\end{equation}
where $\Phi _R$ is an element of the domain of the Dirichlet Laplacian on the cube 
$\Lambda _R\setminus \bigcup _{p\in R_1\Z ^d\cap \Lambda _{3R}}(p+K)$ 
such that $\| \Phi _R\| _2=1$ (cf.\ Theorem (5.25) in \cite{PF92}). 
We take $\Phi _R$ as $\phi _R\psi _R/\| \phi _R\psi _R\|_2$, where $\psi _R$ is the nonnegative and normalized ground state of the Dirichlet Laplacian on $\Lambda _R$ and
 \begin{equation}
  \begin{split}
   \phi _R(x)=
   \left\{
   \begin{array}{lr}
    \Big( 2d_{\infty}\Big( x,\sum_{p\in R^{\mu}\Z ^d\cap \Lambda _R}\Lambda _{R^{\nu}}(p)\Big) R^{-\nu}\Big) \wedge 1 &(d\ge 3),\\[8pt]
    \frac{\ds \Big( \log d_{\infty}(x,\Lambda _R\cap \frac{R\Z ^2}{\sqrt{\log R}})-\frac4{\alpha}\log R\Big) _+}{\ds \log \frac{R}{2\sqrt{\log R}}-\frac4{\alpha}\log R} &(d = 2). 
   \end{array}\right.
  \end{split} \label{testfct}
 \end{equation}
In \eqref{testfct}, $d_{\infty}(\cdot ,\cdot )$ is the distance function with respect to the maximal norm,
$\nu =2/(\alpha -d)$, and $(\cdot )_+$ is the positive part.
Then it is not difficult to see $\| \nabla \Phi _R\| _2^2\le c_4R^{-2}$.
On the event \eqref{lower-event}, we have in addition that
\begin{equation}
 \sum_{q\in \Z ^d}u(x-q-\xi _q)\le \frac{c_5R_1^{d}}{d(x,\sum_{p\in R_1\Z ^d\cap \Lambda _{2R}}\Lambda _1(p))^{\alpha}}
 +c_6R_1^{-(\alpha -d)}\label{ptl-bd}
\end{equation}
in $\Lambda_R$. 
Hence we have 
$$\Big( \Phi _R,\sum_{q\in \Z ^d}u(x-q-\xi _q)\Phi _R\Big) \le c_7R^{-2}.$$
On the other hand, the probability of the event \eqref{lower-event} can be estimated as
 \begin{equation*}
  \begin{split}
   & \log \P(\text{ the event \eqref{lower-event} occurs }) \\
   \ge &\,-\# (R_1\Z ^d\cap \Lambda _{3R})\sum_{q\in \Z ^d\cap \Lambda _{R_1}}\log \P(\xi _0\in \Lambda _1(q)) \\
   & +\sum_{q\in \Z ^d\setminus\Lambda _{2R}}\log (1-\P(|\xi _0|\ge |q|/4))\\
   \ge &\,-c_8R^{d}R_1^{\theta}
  \end{split}
 \end{equation*}
by using $\log (1-X)\ge -2X$ for $0\le X\le 1/2$ in the last line.
Therefore, we have
$$N(c_9R^{-2}) \ge R^{-d} \exp\Bigl( -c_{10} R^dR_1^{\theta} \Bigr) $$
and the proof is finished.
\qed

\begin{rem}
For the manner of taking the function $\phi _R$ in \eqref{testfct} and the event in \eqref{lower-event}, 
we refer the reader to the notion of the ``constant capacity regime'' (cf. Section 3.2.B of \cite{Szn98}).
The same technique is used in Appendix B of \cite{Fuk09a}. 
\end{rem}
\subsection{Compact case}\label{compact}
In this subsection, we adapt the methods in the preceding sections to give a simple proof of the following 
results in \cite{Fuk09a}:

\begin{thm}\label{cpt}
Assume $\Lambda _{r_1}\subset \supp (u)\cup K\subset \Lambda _{r_2}$ for some $0<r_1\le r_2<\infty$ instead of \eqref{single-decay}.
Then we have 
\begin{eqnarray*}
 \log N(\lambda ) \left\{ 
  \begin{array}{ll}
   \sim -(\pi ^2h/\lambda )^{(1+\theta )/2}
   (1+\theta )^{-1}2^{-\theta} & (d=1), \\[8pt]
   \asymp -\lambda ^{-1-\theta /2}(\log (1/\lambda ))^{-\theta /2} & (d=2), \\[5pt]
   \asymp -\lambda ^{-(d/2+\theta /d)} & (d\ge 3)
 \end{array} \right.
\end{eqnarray*}
as $\lambda \downarrow 0$, where $f(\lambda ) \sim g(\lambda )$ means $\lim_{\lambda \downarrow 0}f(\lambda )/g(\lambda )=1$
and $f(\lambda ) \asymp g(\lambda )$ means $0<\varliminf _{\lambda \downarrow 0}f(\lambda )/g(\lambda )\le \varlimsup _{\lambda \downarrow 0}f(\lambda )/g(\lambda )<\infty$.
\end{thm}

\begin{rem}
The assumption on $u$ in this theorem is only for giving a simple proof in the multidimensional case.
If $d=1$, then the assumption in Proposition \ref{1D-upper} is sufficient.
If $d\ge 3$, then this theorem can be extended to the case that the scattering length of $u$ is positive. 
\end{rem}

The proof for $d=1$ is given in Subsection 3.1.
The lower estimate for $d=2$ is given in Subsection 3.3.
To prove the lower estimate for $d\ge 3$, we replace $R^{\nu}$ by $2r_2+1$ in the proof of Proposition \ref{lifshitz-lower}.
Then the rest of the proof is simpler than that of the proposition since
$$\Big( \Phi _R,\sum_{q\in \Z ^d}u(x-q-\xi _q)\Phi _R\Big) =0$$
under the event in \eqref{lower-event} with $R_1=R^{2/d}$.
To prove the upper estimate for $d\ge 3$, we have only to apply the following instead of 
Proposition \ref{ev-bdd} in the proof of Proposition \ref{lifshitz-upper}:

\begin{prop}\label{ev-bdd-cpt}
For sufficiently small $\e _1, \e _2>0$,
there exists a finite constant $c$ such that $\# \{ q\in \Z ^d\cap \Lambda _R : |\xi _q|\ge \e _1R^{2/d}\} \le \e _2R^d$ implies
 \begin{equation}
  \lambda _1\Big( \Big( -\Delta +c_0\sum_{q\in \Z ^d\cap \Lambda _R}1_{B(q+\xi _q,r_0)}\Big) _R^N\Big)
  \ge c/R^2,
 \end{equation}
where $c_0$ and $r_0$ are arbitrarily fixed positive constants.
\end{prop}

\noindent {\it Proof}.
We use the classification
$${\mathcal F}_0= \{ a\in \Lambda _R\cap R^{2/d}\Z ^d 
: \Lambda _{R^{2/d}}(a)\cap \{ q+\xi _q : q\in \Z ^d\cap \Lambda _R\} =\emptyset \}$$
and
$${\mathcal N}_0= \{ a\in \Lambda _R\cap R^{2/d}\Z ^d 
: \Lambda _{R^{2/d}}(a)\cap \{ q+\xi _q : q\in \Z ^d\cap \Lambda _R\} \ne \emptyset \} ,$$
instead of ${\mathcal F}$ and ${\mathcal N}$ in the proof of Proposition \ref{ev-bdd}.
Then we complete the proof by Lemmas \ref{unif-pt-1} and \ref{unif-ef-bdd} 
without using Lemma \ref{unif-pt-2}.
\qed

To prove the upper estimate for $d=2$, we have only to apply the following instead of Proposition \ref{ev-bdd} in the proof of Proposition \ref{lifshitz-upper}:

\begin{prop}\label{ev-bdd-cpt-2d}
For sufficiently small $\e _1, \e _2>0$,
there exists a finite constant $c$ such that $\# \{ q\in \Z ^2\cap \Lambda _R : |\xi _q|\ge \e _1R/\sqrt{\log R}\} \le \e _2R^2$ implies
 \begin{equation}
  \lambda _1\Big( \Big( -\Delta +c_0\sum_{q\in \Z ^2\cap \Lambda _R}1_{B(q+\xi _q,r_0)}\Big) _R^N\Big)
  \ge c/R^2.
 \end{equation}
\end{prop}

To prove this, we replace $R^{2/d}$ by $R/\sqrt{\log R}$ in the proof of Proposition \ref{ev-bdd-cpt}
and we further need to extend Lemma \ref{unif-pt-1} to the 2-dimensional case.
By a simple modification of the proof of Lemma \ref{unif-pt-1}, we have the following, which is 
sufficient for our purpose:

\begin{lem}\label{unif-pt-1-2d} If $d=2$, then we have
 $\inf \{ \lambda _1((-\Delta +c_01_{B(b,r_0)})_R^N) : b\in \Lambda _R\} \ge c/(R^2\log R)$. 
\end{lem}

\section{Critical case} \label{sec-critical}
In this section we discuss the case of $\alpha =d+2$.
By modifying our proof of Proposition \ref{pastur-lower}, we can prove the following:

\begin{prop}\label{critical-lower}
If $\alpha = d+2$, then we have
 \begin{equation}
  \mathop{\varliminf}_{t\uparrow \infty} \frac{\log \NT (t)}{t^{(d+\theta )/(d+2+\theta )}}\ge -K_0(h,C_0),
 \end{equation}
where
  \begin{equation}
   \begin{split}
   & K_0(h,C_0) \\
   =& \inf \Big\{ h\| \nabla \psi \| _2 ^2+\int_{\R ^d}dq \inf_{y\not\in \supp (\psi )-q}
  \Big(   \int_{\R ^d}\frac{dxC_0\psi (x)^2}{|x-q-y|^{d+2}}+|y|^{\theta}\Big) \\
   & \hspace{5mm} :\psi \in W_2^1(\R ^d), \| \psi \| _2=1\Big\} 
   \end{split} \label{critical-bdd}
 \end{equation}
and $W_2^1(\R ^d)=\{ \psi \in L^2(\R ^d) : \nabla \psi \in L^2(\R ^d)\}$.
\end{prop}

\noindent {\it Proof}.
In \eqref{lower-1}, we replace $\psi _R$ by an arbitrary function $\varphi \in H_0^1(\Lambda _R)$ with $\| \varphi \| _2=1$, where $H_0^1(\Lambda _R)$ is the completion of $C_0^{\infty}(\Lambda _R)$ in $W_2^1(\R ^d)$.
Then \eqref{lower-2} is modified as
 \begin{equation*}
  \begin{split}
   & \int _{|q|\le t^{\beta}}\NT _2(t,q)dq \\
   \ge & \mathop{\int}_{|q|\le t^{\beta}}dq\log \mathop{\int}_{y\in [\supp (\varphi ) : R_1+\sqrt{d}/2]^c-q}\frac{dy}{Z(d,\theta )}\\
   & \times \exp \Big( -\int \frac{dx \varphi (x)^2tC_0(1+\e _1)}{\inf \{ |x-q-z-y|^{d+2} : z\in \Lambda _1\}}-|y|^{\theta}\Big) ,
  \end{split} 
 \end{equation*}
where $[A:r]=\{ x\in \R ^d : d(x,A)<r\}$ for any $A\subset \R ^d$ and $r>0$.
We take  $\eta$ as $1/(d+2+\theta )$.
Then, by changing the variables, we see that the right hand side equals
$$t^{d\eta }\mathop{\int}_{|q|\le t^{\beta -\eta }}dq\log \mathop{\int}_{y\in [\supp (\varphi_{\eta} )  : (R_1+\sqrt{d}/2)/t^{\eta}]^c-q}\frac{dy t^{d\eta}}{Z(d,\theta )}\exp (-t^{\theta \eta }\NT _3(y,q;\varphi _{\eta})),$$
where 
$$\NT _3(y,q;\varphi _{\eta})= \int \frac{dx \varphi _{\eta}(x)^2C_0(1+\e _1)}{\inf \{ |x-q-z-y|^{d+2} : z\in \Lambda _{t^{-\eta}}\}}+|y|^{\theta}$$
and $\varphi _{\eta}(x)=t^{d\eta /2}\varphi (t^{\eta}x)$.
We take $R$ as the integer part of ${\mathcal R}t^{\eta}$ for a positive number ${\mathcal R}$
and take $\varphi$ so that $\varphi _{\eta}=\psi$ is a $t$-independent element of $H_0^1(\Lambda _{\mathcal R})$.
Since $t\| \nabla \varphi \| _2^2=t^{(d+\theta )\eta}\| \nabla \psi \| _2^2$ is not negligible, \eqref{lower-4} is modified as
$$\mathop{\varliminf}_{t\uparrow \infty} \frac{\log \NT (t)}{t^{(d+\theta )\eta }}
\ge -h\| \nabla \psi \| _2^2-\mathop{\varlimsup}_{t\uparrow \infty} \int_{|q|\le t^{\beta -\eta }}dq \NT _4(q,t),$$
where
$$\NT _4(q,t) = \inf \Big\{ \mathop{\sup}_{y\in B(y_0,t^{-\gamma})}\NT _3(y,q;\psi ) 
: y_0 \in \Big[ \supp (\psi ) : \frac{R_1+\sqrt{d}/2}{t^{\eta}}+\frac1{t^{\gamma}}\Big] ^c-q\Big\} .$$
Since
$$\mathop{\varlimsup}_{t\uparrow \infty} \NT _4(q,t)
\le \inf_{y\in (\supp (\psi ) )^c-q}\Big( \int \frac{dx \psi (x)^2C_0(1+\e _1)}
{|x-q-y|^{d+2}}+|y|^{\theta}\Big) ,$$
we obtain
$$\mathop{\varliminf}_{t\uparrow \infty} \frac{\log \NT (t)}{t^{(d+\theta )\eta }}
\ge -h\| \nabla \psi \| _2^2-\int_{\R ^d}dy \inf_{y\in (\supp (\psi ))^c-q}\Big( \int \frac{dx \psi (x)^2C_0(1+\e _1)}
{|x-q-y|^{d+2}}+|y|^{\theta}\Big)$$
by the Lebesgue convergence theorem. 
By taking the supremum with respect to $\e _1$, $\psi$ and ${\mathcal R}$, we obtain the result.
\qed

If we apply Donsker and Varadhan's large deviation theory without caring about the topological problems, 
then the formal upper estimate
 \begin{equation}
  \mathop{\varlimsup}_{t\uparrow \infty} \frac{\log \NT (t)}{t^{(d+\theta )/(d+2+\theta )}}\le -K(h,C_0) \label{conj}
 \end{equation}
is expected, where $K(h,C_0)$ is the quantity obtained by removing the restriction $y\not\in \supp (\psi )-q$ 
in the definition \eqref{critical-bdd} of $K_0(h,C_0)$.
For the corresponding Poisson case, this is rigorously established in \^Okura \cite{Oku81}.
In that case, the space $\R ^d$ can be replaced by a $d$-dimensional torus and 
the Feynman-Kac functional becomes a lower semicontinuous functional, 
so that Donsker and Varadhan's theory applies.
However, verifications of both the replacement of the space and the continuity of the functional 
seem to be difficult in our case.

From the conjecture \eqref{conj}, we expect that the quantum effect appears in the leading term.
By Proposition \ref{lifshitz-upper} in Section 3, we can justify this if $d\ge 3$ and $h$ is large:

\begin{prop}\label{critical-3D}
If $d \ge 3$ and $\alpha = d+2$, then we have 
 \begin{equation}
   \mathop{\varlimsup}_{h\to \infty}\mathop{\varlimsup}_{\lambda \to 0}
   \lambda^{(d+\theta )/2}\log N(\lambda) =-\infty .
 \end{equation}
\end{prop}

In the one-dimensional case we can show the same statement with a more explicit bound
$$\mathop{\varlimsup}_{\lambda \to 0}\lambda ^{(1+\theta )/2}\log N(\lambda) \le
   -\frac{\pi ^{1+\theta}h^{(1+\theta )/2}}{(1+\theta )2^{\theta}}$$
by Theorem \ref{lifshitz}, since the leading order does not depend on $\alpha \ge 3$.
In the two-dimensional case we have no such results.
\section{Proof of Theorem \ref{negative}} \label{sec-negative}
\subsection{Upper estimate}
Let $\NT ^-(t)$ be the Laplace-Stieltjes transform of the integrated density of states $N^-(\lambda )$:
$$\NT ^-(t)=\int_{-\infty}^{\infty}e^{-t\lambda}dN^-(\lambda ).$$

To prove the upper estimate, we have only to show the following:

\begin{prop}\label{negative-upper}
Under the condition that $u\ge 0, \sup u=u(0)<\infty$ and $\sup |x|^{\alpha}u(x)<\infty$ for some $\alpha >d$, 
we have
 \begin{equation}
  \mathop{\varlimsup}_{t\uparrow \infty} \frac{\log \NT ^-(t)}{t^{1+d/\theta}}\le u(0)^{1+d/\theta}\int_{|q|\le 1}dq (1-|q|^{\theta}) .
 \end{equation}
\end{prop}

\noindent {\it Proof}.
We use the bound
$$\NT ^-(t)\le \NT _1^-(t)(4\pi th )^{-d/2}$$
as in \eqref{upper-1}, where
$$\NT _1^-(t) = \int_{\Lambda _1}dx\E\Bigg[ \exp \Big( t\sum_{q \in \Z^d} u(\,x-q-\xi_q)\Big) \Bigg] .$$
Here we have used the path integral expression of $\NT ^-(t)$ in Theorem VI.1.1 of \cite{CL90}. 
The assumption required in that theorem will be checked in Lemma \ref{integrable} in Section \ref{Appendix}. 
By replacing the summation by integration, we have
$$\log \NT _1^-(t)\le \int _{\R ^d}dq\log \NT _2^-(t,q),$$
where
$$\NT _2^-(t,q) = \E\Big[ \exp \Big( t \sup_{x\in \Lambda _2}u(x-q-\xi _0)\Big) \Big] .$$
Now we fix an arbitrary small number $\e>0$ and let $C = \sup |x|^{\alpha}u(x)$. 
When $|q| > (1+\e)(u(0)t)^{1/\theta}$, we estimate as 
 \begin{equation}
   \begin{split}
  \NT _2^-(t,q)\le & \exp ( t\sup \{ u(x-y) : x\in \Lambda _2, |y|\ge \delta|q|\} )\\
  & +\exp (tu(0))\P(|\xi _0|\ge (1-\delta)|q|),
   \end{split}\label{nup-1}
 \end{equation}
where $\delta>0$ is chosen to satisfy $(1-\delta)^{\theta+2}(1+\e)^{\theta}=1$. 
For the first term in the right hand side, we use an obvious bound 
$$\sup \{ u(x-y) : x\in \Lambda _2, |y| \ge \delta|q|\} \le C(\delta |q|-\sqrt{d})^{-\alpha}.$$
For the second term, it is easy to see 
$$\P(|\xi _q|\ge (1-\delta)|q|) \le M(\delta, \theta)\exp (-(1-\delta)^{\theta+1}|q|^{\theta})$$
for some large $M(\delta, \theta)>0$. 
Moreover, we have  
$$(1-\delta)^{\theta+1}|q|^{\theta} =\, (1-\delta)^{\theta+2}|q|^{\theta}+ \delta(1-\delta)^{\theta+1}|q|^{\theta}
\ge \, u(0)t + \delta(1-\delta)^{\theta+1}|q|^{\theta}$$
thanks to $|q|>(1+\e)(u(0)t)^{1/\theta}$ and our choice of $\delta$. 
Combining above three estimates, we get 
 \begin{equation}
  \NT _2^-(t,q)\le \exp ( t C(\delta|q|-\sqrt{d})^{-\alpha} )
  (1+M(\delta, \theta)\exp(-\delta(1-\delta)^{\theta+1}|q|^{\theta})) 
 \end{equation}
and thus
 \begin{equation}
  \log \NT _2^-(t,q)\le t C(\delta|q|-\sqrt{d})^{-\alpha}
  +M(\delta, \theta)\exp(-\delta(1-\delta)^{\theta+1}|q|^{\theta}), 
 \end{equation}
using $\log(1+X) \le X $. Since the integral of the right hand side over $\{|q|>(1+\e)(u(0)t)^{1/\theta}\}$ 
is easily seen to be $o(t^{1+d/\theta})$, we can neglect this region. 

For $q$ with $|q|\le (1+\e)(u(0)t)^{1/\theta}$, we estimate as
 \begin{equation}
  \begin{split}
  \NT _2^-(t,q)\le &\exp ( t\sup \{ u(x-y) : x\in \Lambda _2, |y|\ge L\} )\\ 
  & +\exp (tu(0))\P(|q+\xi _0|\le L),
 \end{split} \label{nup-2}
 \end{equation}
where $L = 2\e(u(0)t)^{1/\theta}$.
We use obvious bounds 
$$\sup \{ u(x-y) : x\in \Lambda _2, |y|\ge L\} \le C(L-\sqrt{d})_+^{-\alpha}$$ 
for the first term and 
$$\P(|q+\xi _0| \le L)\le \exp (-(|q|-L)_+^{\theta})|B(0,L)|/Z(d,\theta )$$
for the second term. 
Note also that we have 
$$tc(L-\sqrt{d})_+^{-\alpha}\le tu(0)-(|q|-L)_+^{\theta}$$
for large $t$, from $|q| \le (1+\e)(u(0)t)^{1/\theta}$ and our choice of $L$. 
Using these estimates, we obtain 
 \begin{equation*}
  \begin{split}
   &\, \int_{|q|\le (1+\e)(u(0)t)^{1/\theta}}dq\log \NT _2^-(t,q)\\
   \le &\, \int_{|q|\le (1+\e)(u(0)t)^{1/\theta}}dq\Big\{ \log \Big( \frac{|B(0,L)|}{Z(d,\theta )}+1\Big)
   +tu(0)-(|q|-L)_+^{\theta}\Big\}.
  \end{split}
 \end{equation*}
By changing the variable and taking the limit, we arrive at 
$$\mathop{\varlimsup}_{t\uparrow \infty} \frac{\log \NT (t)}{t^{1+d/\theta}}
\le u(0)^{1+d/\theta}\int_{|q|\le 1+\e}dq\{ 1-(|q|-2\e)_+^{\theta}\} .$$
This completes the proof of Proposition \ref{negative-upper} since $\e > 0$ is arbitrary.
\qed

\subsection{Lower estimate}
To prove the lower estimate, we have only to show the following:

\begin{prop}\label{negative-lower}
Suppose $u\ge 0, \sup u=u(0)<\infty$, and $u(x)$ is lower semicontinuous at $x=0$. Then we have
 \begin{equation}
  \mathop{\varliminf}_{t\uparrow \infty} \frac{\log \NT ^-(t)}{t^{1+d/\theta}}\ge u(0)^{1+d/\theta}\int_{|q|\le 1}dq (1-|q|^{\theta}) .
 \end{equation}
\end{prop}

\noindent {\it Proof}.
For any $\e >0$, there exists $R_{\e}>0$ 
such that 
\begin{equation}
	u(x)\ge u(0)-\e \textrm{ for } |x|<R_{\e}\label{lsc}
\end{equation}
by the lower semicontinuity of $u$.
We use the bound
$$\NT ^-(t)\ge \exp (-th\| \nabla \psi _{\e}\| _2)\NT _1^-(t),$$
for any $\psi _{\e}\in C_0^{\infty} (\Lambda _{\e})$ such that the $L^2$-norm of $\psi _{\e}$ is 1, where
\begin{equation}
\NT _1^-(t) = \E\Bigg[\exp \Big( t\sum_{q \in \Z^d} \inf_{x\in \Lambda _{\e}}u(\,x-q-\xi_q)\Big) \Bigg] .\label{ntm1}
\end{equation}
This is proven by the same estimate as used in \eqref{lower-1}.
We take $\psi _{\e}$ as the nonnegative and normalized ground state of the Dirichlet Laplacian on the cube $\Lambda _{\e}$.
Since a sufficient condition for $\sup_{x\in \Lambda _{\e}}|x-q-\xi _q|\le R_{\e}$ is $|q+\xi _q|\le R_{\e}-\e \sqrt{d}/2$, 
we restrict the expectation to this event and deduce  from~\eqref{lsc} that
$$\log \NT _1^-(t)\ge \sum_{q\in \Z ^d}\log \int_{|q+y|\le R_{\e}-\e \sqrt{d}/2}\frac{dy}{Z(d,\theta )}\exp (t(u(0)-\e )-|y|^{\theta}).$$
Since a sufficient condition for $\inf \{ u(0)-\e )-|y|^{\theta}\le R_{\e} : |q+y|\le R_{\e}-\e \sqrt{d}/2\} \ge 0$ is 
$|q|\le \{ t(u(0)-\e )\} ^{1/\theta}-R_{\e}+\e \sqrt{d}/2$, we restrict the range of $q$ and deduce 
\begin{equation*}
  \begin{split}
   & \log \NT _1^-(t)\\
   & \ge \int_{|q|\le h(t)}\Big\{ c'\log \frac{|B(0,R_{\e}-\e \sqrt{d}/2)|}{Z(d,\theta )}+t(u(0)-\e )-(|q|+R_{\e}-c))^{\theta}\Big\} \\
   & =h(t)^d\int_{|q|\le 1}\Big\{ c'\log \frac{|B(0,R_{\e}-\e \sqrt{d}/2)|}{Z(d,\theta )}+t(u(0)-\e )-(h(t)|q|+R_{\e}+c))^{\theta}\Big\} 
  \end{split}
 \end{equation*}
for large $t$ and small $\e$, where $h(t)=\{ t(u(0)-\e )\} ^{1/\theta}-R_{\e}-c$ and $c$ and $c'$ are positive constants.
Then we obtain
$$\mathop{\varliminf}_{t\uparrow \infty} \frac{\log \NT ^-(t)}{t^{1+d/\theta}}\ge (u(0)-\e )^{1+d/\theta}
\int_{|q|\le 1}dq (1-|q|^{\theta}).$$
Since $\e$ is arbitrary,
this completes the proof of Proposition \ref{negative-lower}.
\qed


\section{Asymptotics for associated Wiener integrals} \label{sec-FK}
In the previous work \cite{Fuk09a}, the asymptotic behaviors of the integrated density of states 
were derived from those of certain Wiener integrals. 
In this section, we recall the connection and derive the asymptotic behaviors of 
the associated Wiener integrals in our settings. 
Let $h=1/2$ for simplicity and $E_{x}$ denote the expectation with respect to the standard 
Brownian motion $(B_s)_{0 \le s \le \infty}$ starting at $x$. 
Then the Laplace-Stieltjes transform of the integrated density of states can be expressed as follows 
(cf. Chapter VI of \cite{CL90}): 
\begin{equation}
 \begin{split}
 \NT(t) = (2\pi t)^{-d/2} & \int_{\Lambda _1}dx\E \otimes E_{x}
 \biggl[ \exp \biggl\{-\int_0^t \sum_{q \in \Z^d}u(B_s-q-\xi_q) ds \biggr\} \\
 & : B_s\not\in \bigcup _{q\in \Z ^d}(q+\xi _q+K) \text{ for }0\le s\le t 
 \biggl| \, B_t = x \biggr].
 \end{split}\label{FK}
\end{equation}
We can also express $\NT^{-}(t)$ in the same form by changing the sign of $u$ and 
setting $K=\emptyset$ in the right hand side. 
In view of \eqref{FK}, $\NT(t)$ seems, and indeed will be proven below, 
to be asymptotically comparable to the Wiener integral 
\begin{equation}
 \begin{split}
 S_{t,\,x} = \E \otimes E_{x} \biggl[ & \exp \biggl\{-\int_0^t \sum_{q \in \Z^d}u(B_s-q-\xi_q) ds \biggr\} \\
 & : B_s\not\in \bigcup _{q\in \Z ^d}(q+\xi _q+K) \text{ for }0\le s\le t \biggr],
 \end{split} \label{sp} 
\end{equation}
which was the main object in \cite{Fuk09a}. 
This quantity is of interest itself since not only it gives the average of the solution of a heat equation 
with random sinks but also can be interpreted as the annealed survival probability of the Brownian motion 
among killing potentials. 
Similarly, $\NT ^-(t)$ is asymptotically comparable to the average of the solution 
\begin{equation}
 S_{t,\,x}^- = \E \otimes E_{x} \biggl[ \exp \biggl\{ \int_0^t \sum_{q \in \Z^d}
 u(B_s-q-\xi_q) ds \biggr\} \biggr], \label{sp-n}
\end{equation}
of a heat equation with random sources which 
can also be interpreted as the average number of the branching Brownian motions 
in random media. 
We refer the readers to \cite{GM90,GK04,Szn98} about the interpretations of $S_{t,\,x}$ and $S_{t,\,x}^-$. 
The connection between the asymptotics of $\NT (t)$ and $S_{t,\,x}$ can be found in the literature 
for the case that $\{ q+\xi _q \} _q$ is replaced by an $\R ^d$-stationary random field 
(see e.g.\ \cite{Nak77}, \cite{Szn90a}).
However our case is only $\Z ^d$-stationary. 

We first prepare a lemma which gives upper bounds on $\log S_{t,\,x}$ and $\log S_{t,\,x}^-$
in terms of $\log \NT (t)$ and $ \log\NT ^-(t)$, respectively. 
We shall state the results only for $x \in \Lambda_1$ since they automatically extend to 
the whole space by the $\Z ^d$-stationarity. 
\begin{lem}\label{sn-bdd}
For any $x \in \Lambda_1$ and $\e > 0$, we have
 \begin{equation}
  \log S_{t,\,x}\le \log \NT (t-\e)(1+o(1)) \label{sn-bdd1}
 \end{equation}
and
 \begin{equation}
  \log S_{t,\,x}^-\le \log \NT ^-(t-t^{-2d/\theta})(1+o(1)) \label{sn-bdd2}
 \end{equation}
as $t \to \infty$. 
\end{lem}

\noindent {\it Proof}.
We give the proof of \eqref{sn-bdd2} first. Let $V_{\xi}(x)$ denotes the potential 
$\sum_{q \in \Z^d} u(x-q-\xi_q)$ for simplicity. 
We divide the expectation as 
\begin{equation}
 \begin{split}
  S_{t,\,x}^- = &\E \otimes E_{x} \biggl[ \exp \biggl\{ \int_0^t V_{\xi}(B_s) ds \biggr\}:\, 
  \sup_{0 \le s \le t}|B_s|_{\infty} < [t^{1+d/\theta}] \biggr]\\
  &+\sum_{n > [t^{1+d/\theta}]}\E \otimes E_{x} \biggl[ \exp \biggl\{ \int_0^t 
  V_{\xi}(B_s) ds \biggr\}:\, n-1 \le \sup_{0 \le s \le t}|B_s|_{\infty} <n \biggr].\label{divided}
 \end{split}
\end{equation}

The summands in the second term can be bounded from above by 
\begin{equation}
 \begin{split}
  &\E\biggl[ \exp\biggl\{t\sup_{y \in \Lambda_{2n}} V_{\xi}(y) \biggr\} \biggr]
  P_x\left(n-1 \le \sup_{0 \le s \le t}|B_s|_{\infty}\right)\\
  \le &c_1n^d\E\biggl[ \exp\biggl\{t\sup_{y \in \Lambda_1} V_{\xi}(y) \biggr\} \biggr]
  \exp\{-c_2n^2/t\}\\
  \le & c_1n^d \exp\{c_3t^{1+d/\theta}-c_2n^2/t\},
 \end{split}
\end{equation}
where we have used a standard Brownian estimate (cf.\ \cite{IM74} Section 1.7) and the $\Z^d$-stationarity in the second line, 
and Lemma \ref{integrable} below in the third line.
Then, it is easy to see that the second term in \eqref{divided} is bounded from above 
by a constant and hence it is negligible compared with $\NT^-(t)$. 

Now let us turn to the estimate of the first term in \eqref{divided}. 
Note first that we can derive an upper large deviation bound 
\begin{equation}
 \P\biggl( \sup_{y \in \Lambda_{[t^{1+d/\theta}]}} V_{\xi}(y) \ge v \biggr)
 \le [t^{1+d/\theta}]^d\P\biggl( \sup_{y \in \Lambda_1} V_{\xi}(y) \ge v \biggr)
 \le \exp (-c_4 v^{1+\theta/d})\label{LDP}
\end{equation}
which is valid for all sufficiently large $t$ and $v \ge t$, 
from the exponential moment estimate in Lemma \ref{integrable} below. 
Using this estimate, we get 
\begin{equation}
 \begin{split}
  &\E \otimes E_{x} \biggl[ \exp \biggl\{ \int_0^t V_{\xi}(B_s) ds \biggr\}:\, 
  \sup_{0 \le s \le t}|B_s|_{\infty} < [t^{1+d/\theta}],\\
  &\qquad\qquad\qquad\qquad\qquad\qquad
  \sup_{y \in \Lambda_{2[t^{1+d/\theta}]}} V_{\xi}(y) \ge t^{2d/\theta} \biggr]\\
  \le &\E\biggl[ \exp\biggl\{t\sup_{y \in \Lambda_{2[t^{1+d/\theta}]}} 
  V_{\xi}(y) \biggr\}:\,
  \sup_{y \in \Lambda_{2[t^{1+d/\theta}]}} V_{\xi}(y) \ge t^{2d/\theta} \biggr]\\
  \le &\sum_{n \ge t^{2d/\theta}} \exp\{tn\}
  \P\biggl(n-1 \le \sup_{y \in \Lambda_{2[t^{1+d/\theta}]}} V_{\xi}(y) < n\biggr)\\  
  \le &\sum_{n \ge t^{2d/\theta}}\exp\set{tn - c_4 (n-1)^{1+\theta/d} }.\\
 \end{split}
\end{equation}
Since the last expression converges to 0 as $t \to \infty$, 
we can restrict ourselves on the event $\{\sup V_{\xi}(x) \le t^{2d/\theta}\}$. 
Hereafter, we let $T=[t^{1+d/\theta}]$ since its exact form will be irrelevant in the sequel. 
Then, the Markov property at time $\e = t^{-2d/\theta}$ yields 
 \begin{equation}
  \begin{split}
   &\E \otimes E_{x} \biggl[ \exp \biggl\{ \int_0^t V_{\xi}(B_s) ds \biggr\}:\, 
   \sup_{0 \le s \le t}|B_s|_{\infty} < T, 
   \sup_{y \in \Lambda_{2T}} V_{\xi}(y) < t^{2d/\theta} \biggr]\\
   \le &e \int _{\Lambda_{2T}}\frac{dy}{(2\pi \e )^{d/2}}\exp \Bigl( -\frac{|x-y|^2}{2\e}\Bigr) \\
   & \times \E \otimes E_{y} \biggl[ \exp \biggl\{\int_0^{t-\e} V_{\xi}(B_s) ds \biggr\} :\, 
   \sup_{0 \le s \le t-\e}|B_s|_{\infty} < T \biggr]\\
   \le &\frac{e}{(2\pi \e )^{d/2}}\int _{\Lambda _{2T}}dy\int _{\Lambda _{2T}}dz \E  
   [ \exp (-(t-\e )H^{-,\,D}_{\xi , \,2T})(y,z) ], \label{sn-bdd-p1}
  \end{split}
 \end{equation}
where $\exp (-tH^{-,\,D}_{\xi , \,2T})(x,y)$, $t>0$,
$x, y\in \Lambda _{2T}$, is the integral kernel of the heat semigroup generated by the self-adjoint operator $H^-_{\xi}$ on the $L^2$-space on the cube $\Lambda _{2T}$ with the Dirichlet boundary condition.

Finally, we use the estimate
$$\exp (-tH^D_{\xi ,\, 2T})(y,z)\le \bigl\{ \exp (-tH^D_{\xi,\,2T})(y,y)\exp (-tH^D_{\xi, \,2T})(z,z)\bigr\}^{1/2}$$
for the kernel of self-adjoint semigroup 
and the Schwarz inequality to dominate the right hand side in \eqref{sn-bdd-p1} 
by $T^{2d}\NT^-(t-\e )$ multiplied by some constant. 

Combining all the estimates above, we finish the proof of \eqref{sn-bdd2}. 
We can also prove \eqref{sn-bdd1} in the same way as \eqref{sn-bdd-p1}. 
However it is much simpler since we do not have to care about $\sup V_{\xi}(\,\cdot\,)$ 
and thus we omit the details. 
\qed

The next lemma gives the converse relation between $\log S_{t,\,x}$ and $\log \NT (t)$, 
while the lower estimate of $\log S_{t,\,x}^-$ will be derived directly. 
(See the proof of Theorem \ref{F-K}.)
\begin{lem}\label{ns-bdd}
For any $x \in \Lambda_1$ and $\e > 0$, we have
\begin{equation}
 \log \NT (t) \le \log S_{t-\e,\,x}^{v,K'}(1+o(1)) \label{ns-bdd1}
\end{equation}
as $t\to \infty$, where $S_{t,\,x}^{v,K'}$ is the expectation defined by replacing $K$ and $u$ by $K'=\{ x\in K : d(x,K^c)\ge \sqrt{d} \}$ and $v(y)=\inf \{ u(y-x+z) : z\in \Lambda _1\}$ 
respectively in \eqref{sp}. 
Note that if $u$ is a function satisfying the conditions in Theorem \ref{pastur} or \ref{lifshitz}, 
then so is $v$.
\end{lem}

\noindent {\it Proof}.
Let $\e > 0$ be an arbitrarily small number. 
By the Chapman-Kolmogorov identity, we have
\begin{equation*}
 \begin{split}
 \NT (t) \le (2\pi \e )^{-d/2}\int _{\Lambda _1}dz
 \E \otimes E_{z} \biggl[ & \exp \biggl\{-\int_0^{t-\e} \sum_{q \in \Z^d}u(B_s-q-\xi_q) ds \biggr\} \\
 & : B_s\not\in \bigcup _{q\in \Z ^d}(q+\xi _q+K) \text{ for }0\le s\le t-\e \biggr].
 \end{split} 
\end{equation*}
The right hand side is dominated by $(2\pi \e )^{-d/2} S_{t-\e,\,x}^{v,K'}$ 
and the proof of \eqref{ns-bdd1} is completed.
\qed

We now state our results on the asymptotics of $S_{t,\,x}$ and $S_{t,\,x}^-$:

\begin{thm}\label{F-K}
{\upshape{(i)}} Assume $d=1$ and \eqref{positive} if $\alpha \le 3$.
Then we have 
\begin{equation}
\begin{split} 
 \log S_{t,\,x} \left\{ 
  \begin{array}{ll}
   \sim -t^{(1+\theta )/(\alpha+\theta )}
   \ds{\int_{\R}dq \inf_{y\in \R}\Big( \frac{C_0}{|q+y|^{\alpha}}+|y|^{\theta}\Big)} & (1 < \alpha < 3), \\[8pt]
   \asymp -t^{(1+\theta )/(3+\theta )} & (\alpha = 3), \\[5pt]
   \ds \sim {-t^{(1+\theta )/(3+\theta )}\frac{3+\theta}{1+\theta}\big( \frac{\pi ^2}8 \big) ^{(1+\theta )/(3+\theta )}}
   & (\alpha >3)
 \end{array} \right. \label{F-K-1}
\end{split}
\end{equation}
as $t \to \infty$, where $f(t) \sim g(t)$ means $\lim_{t\to \infty}f(t)/g(t)=1$
and $f(t) \asymp g(t)$ means $0<\varliminf _{t\to \infty}f(t)/g(t)\le \varlimsup _{t\to \infty}f(t)/g(t)<\infty$.

{\upshape{(ii)}} Assume $d=2$ and \eqref{positive} if $\alpha \le 4$.
Then we have  
\begin{equation}
\begin{split} 
\log S_{t,\,x} \left\{ 
  \begin{array}{ll}
   \sim -t^{(2+\theta )/(\alpha+\theta )}
   \ds{\int_{\R ^2}dq \inf_{y\in \R ^2}\Big( \frac{C_0}{|q+y|^{\alpha}}+|y|^{\theta}\Big)} & (2 < \alpha < 4), \\[8pt]
   \asymp -t^{(2+\theta )/(4+\theta )} & (\alpha = 4), \\[5pt]
   \asymp {-t^{( 2+\theta )/(4+\theta )} (\log t)^{-\theta /(4+\theta )}} & (\alpha >4)
 \end{array} \right. \label{F-K-2}
\end{split}
\end{equation}
as $t \to \infty$. 

{\upshape{(iii)}} Assume $d \ge 3$ and \eqref{positive} if $\alpha \le d+2$.
Then we have 
\begin{equation}
\begin{split} 
 \log S_{t,\,x} \left\{ 
  \begin{array}{ll}
   \sim -t^{(d+\theta )/(\alpha+\theta )}
   \ds{\int_{\R ^d}dq \inf_{y\in \R ^d}\Big( \frac{C_0}{|q+y|^{\alpha}}+|y|^{\theta}\Big)} & (d < \alpha < d+2), \\[10pt]
   \asymp {-t^{(d+\theta\mu )/(d+2+\theta\mu )}} & (\alpha \ge d+2)
 \end{array} \right. \label{F-K-3}
\end{split}
\end{equation}
as $t \to \infty$, where $\mu =2(\alpha -2)/(d(\alpha -d))$ as in Theorem \ref{lifshitz}. 

{\upshape{(iv)}} Assume $\sup u=u(0)<\infty$ and the existence of $R_{\e}>0$ for any $\e >0$
such that $\essinf _{B(R_{\e})}u\ge u(0)-\e$. 
Then we have 
\begin{equation}
 \log S_{t,\,x}^{-} \sim t^{1+d/\theta} u(0)^{1+d/\theta}\int_{|q|\le 1}dq (1-|q|^{\theta}) \label{F-K-4}
\end{equation}
as $t \to \infty$. 
\end{thm}

\noindent {\it Proof}.
We first consider the corresponding results for $\NT (t)$ and $\NT ^-(t)$: the estimates \eqref{F-K-1}--\eqref{F-K-4} 
with $S_{t,\,x}$ and $S_{t,\,x}^{-}$ replaced by $\NT (t)$ and $\NT ^-(t)$, respectively.  
These are already proven in earlier sections except for the case of $\alpha > d+2$ and $d\ge 2$. 
The results for the remaining case follow from Propositions \ref{lifshitz-upper} 
and \ref{lifshitz-lower} and Abelian theorems in~\cite{Kas78}. 
Then by Lemma \ref{sn-bdd}, we obtain the upper estimates of $S_{t,\,x}$ and $S_{t,\,x}^{-}$.
For the lower estimates of $S_{t,\,x}$, we set $u^{\#}(y)=\sup \{ u(y+x+z) : z\in \Lambda _1\} 1_{B(R_1)^c}(y)+1_{B(R_1)}(y)$  with $R_1\ge 0$.
If $u$ satisfies the conditions in Theorems \ref{pastur} and \ref{lifshitz}, and $R_1$ is sufficiently large, then $u^{\#}$ also satisfies the same conditions.
Therefore we obtain the corresponding lower estimates of $\NT (t)$
where $K$ is replaced by $B(R_2)$ with any $R_2\ge R_1$ and $u$ is replaced by $u^{\#}$.
Then by Lemma \ref{ns-bdd}, we obtain the corresponding lower estimates of $S_{t,\,x}^{v^{\#},B(R_2+\sqrt{d})}$, where $v^{\#}(y)=\inf \{ u^{\#}(y-x+z) : z\in \Lambda _1\}$.
Since $K\subset B(R_2+\sqrt{d})$ and $v^{\#}\ge u$ on $B(R_2)^c$ for some $R_2\ge R_1$, we obtain the corresponding lower estimates of $S_{t,\,x}$.
For the lower estimate of $S_{t,\,x}^-$, we restrict the expectation to the event $B_s\in \Lambda _{\e}$ 
for any $s\in [1,t]$ to obtain
$$ S_{t,\,x}^- \ge \int_{\Lambda _{\e}}dy e^{\Delta /2}(x,y)
\int_{\Lambda _{\e}}dz e^{(t-1)\Delta ^D_{\e}/2}(y,z)\NT _1^-(t-1)\ge c_1e^{-c_2t}\NT _1^-(t-1), $$
where $\NT _1^-(t)$ is the function defined in \eqref{ntm1}, and 
$\exp (t\Delta /2)(x,y)$, $(t,x,y)\in (0,\infty )\times \R ^d\times \R ^d$ and $\exp (t\Delta ^D_{\e}/2)(x,y)$, 
$(t,x,y)\in (0,\infty )\times \Lambda _{\e}\times \Lambda _{\e}$ are the 
integral kernels of the heat semigroups 
generated by the Laplacian and the Dirichlet Laplacian on $\Lambda _{\e}$, respectively, multiplied by $-1/2$.
Therefore the lower estimate of $S_{t,\,x}^-$ is given by our proof of Proposition \ref{negative-lower}.
\qed


\section{Appendix}\label{Appendix}
We here state and prove two lemmas which we used before. 
The first one is to define the integrated density of states $N(\lambda)$ and 
to represent it by the Feynman-Kac formula: 

\begin{lem} \label{kato}
 Let $u$ be a nonnegative function belonging to the class $K_d$ and satisfying \eqref{single-decay}.
 Let $\xi =(\xi_q)_{q \in \Z^d}$ be a collection of independently and identically
distributed $\R^d$-valued random variables satisfying \eqref{single-prob}.
 Then almost all sample functions of the random field defined by $V_{\xi}(x)=\sum_{q \in \Z^d} u(x-q-\xi_q)$ belong to the class $K_{d,loc}$. 
\end{lem}

\noindent {\it Proof}.
For any $\e , \delta >0$, by the Chebyshev inequality, we have
$$\P (|\xi _q|\ge |q|^{\e})\le \E [(|\xi _q|/|q|^{\e})^{\delta}]\le c_1/|q|^{\e \delta}.$$
For any $\e$, there exists $\delta$ such that
$$\sum_{q\in {\Z}^d}\P (|\xi _q|\ge |q|^{\e})<\infty .$$
By the Borel-Cantelli lemma, for almost all $\xi$, we have $N_{\xi}\in {\N}$ such that $|\xi _q|<|q|^{\e} < |q|/3$ for any $q\in \Z -B(N_{\xi})$.
By the condition \eqref{single-decay} we also have $R_{\e}$ such that $u(x)\le (C_0+\e )/|x|^{\alpha}$ for any $x\in B(R_{\e})^c$.
We now take $R>0$ arbitrarily.
If $x\in B(R)$ and $q\in \Z ^d-B(3(R\vee R_{\e})\vee N_{\xi})$, then 
$$|x-q-\xi _q|\ge |q|-|\xi _q|-|x|\ge |q|/3\ge R_{\e}$$
and 
$$V_{\xi}(x)\le \sum_{q\in \Z ^d\cap B(3(R\vee R_{\e})\vee N_{\xi})}u(x-q-\xi_q)+c_2.$$
Since the right hand side is a finite sum, we have $1_{B(R)}V_{\xi}\in K_d$.
Since $R$ is arbitrary, we complete the proof.
\qed

The second is to define the integrated density of states $N^-(\lambda )$ and represent it by the 
Feynman-Kac formula. The following is enough to apply Theorem VI.1.1 in \cite{CL90}. 
This lemma was also used in~\eqref{LDP}. 
\begin{lem} \label{integrable}
 Let $u$ be a bounded nonnegative function satisfying \eqref{single-decay}.
 Then there exist finite constants $c_1$ and $c_2$ such that 
$$\E \Big[ \exp \Big( r\sup_{x\in \Lambda _1}V_{\xi}(x)\Big) \Big] \le c_1\exp (c_2r^{1+d/\theta})$$
for any $r\ge 0$,
where $\xi$ and $V_{\xi}$ are same as in the last lemma. 
\end{lem}

\noindent {\it Proof}.
We first dominate as
$$\log \E \Big[ \exp \Big( r\sup_{x\in \Lambda _1}V_{\xi}(x)\Big) \Big] 
\le \int_{\R ^d}\log I(q)dq,$$
where
$$I(q)=\E \Big[ \exp \Big( r\sup_{x\in \Lambda _2}u(x-q-\xi _0)\Big) \Big] .$$
For sufficiently large $R>0$, we have $u(x)\le 2C_0|x|^{-\alpha}$ for $|x|\ge R_0$.
A sufficient condition for $\inf_{x\in \Lambda _2}|x-q-\xi _0|\ge R$ is $|q+\xi _0|\ge R+\sqrt{d}$.
Then, for $q\in B(2(R+\sqrt{d}))^c$, we dominate as
\begin{equation*}
 \begin{split}I(q)\le & \E \Big[ \exp \Big( \sup_{x\in \Lambda _2}\frac{2rC_0}{|x-q-\xi _0|^{\alpha}}\Big) :  |q+\xi _0|\ge \frac{|q|}2 \Big] 
+e^{r\sup u} \P \Big( |q+\xi _0|<\frac{|q|}2 \Big) \\
  \le & \exp \Big( \frac{2rC_0}{(|q|/2-\sqrt{d})^{\alpha}}\Big) (1 + c_1\exp (r\sup u-c_2|q|^{\theta})) 
 \end{split}
\end{equation*}
Since $\log (1+X)\le X$ for any $X\ge 0$, we have
\begin{equation*}
 \begin{split}& \int_{B(2(R+\sqrt{d}))^c}\log I(q)dq \\
\le & \int_{B(2(R+\sqrt{d}))^c}\frac{2rC_0}{(|q|/2-\sqrt{d})^{\alpha}}dq+\int_{B(2(R+\sqrt{d}))^c}c_1\exp (r\sup u-c_2|q|^{\theta}))dq\\
 \le & \frac{c_3r}{R^{\alpha -d}}+c_4\exp (r\sup u-c_5R^{\theta}).
 \end{split}
\end{equation*}
By a simple uniform estimate, we have
$$\int_{B(2(R+\sqrt{d}))}\log I(q)dq\le c_6r\sup uR^d.$$
Setting $R=(r\sup u/c_5)^{1/\theta}$, 
we have
$$\int \log I(q)dq\le c_7r^{1+d/\theta}$$
 for sufficiently large $r>0$.
\qed



\end{document}